%%%%%%%%%%%%%%% PLAIN TeX %%%%%%%%%%%%%%%

%-----------------HKLS29.4.tex ---------
%------------22/11/2004/Tokyo-----------

%-----Paper Format---------
\hsize=13.5cm
\vsize=20.3cm
\hoffset=1.2cm
\voffset=2.5cm
\baselineskip=13pt

\nopagenumbers

%-----AMS Symbols-------
\input amssym.def
\input amssym.tex

%-----Special Fonts-----
\font\srm=cmr8 

\font\title=cmr12 at 15pt
\def\L{{\font\lmat=cmmi12 at 15pt\lmat L}}
\def\sL{{\font\slmat=cmmi8 \slmat L}}

\font\teneusm=eusm10
\font\seveneusm=eusm7
\font\fiveeusm=eusm5
\newfam\eusmfam
\def\eusm{\fam\eusmfam\teneusm}
\textfont\eusmfam=\teneusm
\scriptfont\eusmfam=\seveneusm
\scriptscriptfont\eusmfam=\fiveeusm

\font\teneufm=eufm10
\font\seveneufm=eufm7
\font\fiveeufm=eufm5
\newfam\eufmfam

\textfont\eufmfam=\teneufm
\scriptfont\eufmfam=\seveneufm
\scriptscriptfont\eufmfam=\fiveeufm

%----Small Macros----
\def\varGamma{{\mit \Gamma}}
\def\txt#1{{\textstyle{#1}}}
\def\scr#1{{\scriptstyle{#1}}}
\def\r#1{{\rm #1}}
\def\B#1{{\Bbb #1}}
\def\e#1{{\eusm #1}}

\def\sgn{{\rm sgn}}

%----Headers-------
\def\rightheadline{\hfil{\srm
UNIFORM BOUND FOR HECKE \sL-FUNCTIONS}
\hfil\tenrm\folio}
\def\leftheadline{\tenrm\folio\hfil{\srm
M. JUTILA and Y. MOTOHASHI}\hfil}
\def\emptyheadline{\hfil}
\headline{\ifnum\pageno=1 \emptyheadline\else
\ifodd\pageno \rightheadline \else \leftheadline\fi\fi}

%-------Text---------
\centerline{\title Uniform bound for Hecke \L-functions}
\vskip 1cm
\centerline{by}
\bigskip
\centerline{MATTI JUTILA and YOICHI MOTOHASHI}
\vskip 2cm
\centerline{\bf 1. Introduction}
\bigskip
\noindent
Our principal aim in the present article is to establish a 
uniform hybrid bound  for individual values on the
critical line of Hecke $L$-functions
associated with cusp forms over the full modular group. This is rendered
in the statement that for $t\ge0$
$$
\eqalignno{
H_j\left(\txt{1\over2}+it\right)&\ll
\left(\kappa_j+t\right)^{1/3+\varepsilon},&(1.1)\cr
H_{j,k}\left(\txt{1\over2}+it\right)&\ll\hskip 2pt 
\left( k+t\right)^{1/3+\varepsilon},&(1.2)\cr
}
$$
with the common notation to be made precise in the course of discussion. 
\medskip
Most of arithmetically significant Dirichlet series such as the Riemann
zeta-function $\zeta(s)$, Dirichlet $L$-functions, and Hecke $L$-functions
associated with various cusp forms satisfy Riemannian functional equations
connecting values at $s=\sigma+it$ and $1-s$ of respective functions. Essentially
best possible estimates for these functions near the line $\sigma=1$ and
$\sigma=0$ can usually be deduced from the definition of respective functions and
their functional equations. From this, bounds in the critical strip $0<\sigma<1$,
in particular on the critical line $\sigma={1\over2}$, follow readily via the
Phragm\'en--Lindel\"of convexity principle; thus they are called convexity
bounds. In general, there is a quantity $B(g,t)$ characterising the size
of a function $g\left({1\over2}+it\right)$ of the above kind in a given $t$-range
in such a way that the convexity bound is stated as
$$
g\left(\txt{1\over2}+it\right)\ll B(g,t)^{1/2+\varepsilon},\quad t> 0,
\eqno(1.3)
$$
with the usual usage of the symbol $\varepsilon$ 
(see Convention 1 at the end of this section).  For instance,
$B(\zeta,t)=t^{1/2}$, or perhaps more naturally $B(\zeta^2,t)=t$. In view of the
generalised Lindel\"of Hypothesis asserting that the exponent on the right of
$(1.3)$ be $\varepsilon$, any improvement upon $(1.3)$, i.e., subconvexity
bounds are of considerable interest. One may call $\varpi<{1\over2}$ a Lindel\"of
constant, provided that $(1.3)$ holds with $\varpi+\varepsilon$ in place of
${1\over2}+\varepsilon$. The classical Lindel\"of constant for $\zeta^2$ is
$\varpi={1\over3}$, which has been successively improved, though not very
drastically. A natural task would then be to achieve at least the same for
wide classes of Dirichlet series. We shall consider this fundamental problem
dealing mainly with Hecke $L$-functions associated with real analytic cusp
forms.
\medskip
To this end, we shall first make our objects precise; for details we refer
to the monograph [23]. Thus, let $\varGamma$ be the full modular group
$\r{PSL}_2(\B{Z})$; throughout the sequel we shall work with $\varGamma$,
although our argument appears to be effective in a considerably 
general setting. Let $L^2\!\left(\varGamma\backslash{\Bbb H}\right)$ 
be the Hilbert space composed of all $\varGamma$-automorphic functions on the
hyperbolic upper half plane $\B{H}=\left\{x+iy\,:\, x\in\B{R},\,y>0\right\}$, 
which are square integrable over the quotient $\varGamma\backslash{\Bbb H}$
with respect to the hyperbolic measure. If a function in
$L^2\!\left(\varGamma\backslash{\Bbb H}\right)$ is an eigenfunction of the
hyperbolic Laplacian $\e{L}=-y^2\left(\partial_x^2 +\partial_y^2\right)$, then it
is called a real analytic cusp form. The subspace spanned by all such functions 
has a maximal orthonormal system $\left\{\psi_j\,:\, j=1,2,\ldots\right\}$, where
$\e{L}\psi_j=\left({1\over4}+\kappa_j^2\right)\psi_j$ with
$0<\kappa_1\le\kappa_2\le\cdots$, and
$$
\psi_j(x+iy)=\sqrt{y}\sum_{\scr{n=-\infty}\atop\scr{n\ne0}}
^\infty\varrho_j(n)K_{i\kappa_j}\left(2\pi|n|y\right)
\exp(2\pi inx),\quad x+iy\in\B{H},\eqno(1.4)
$$
with $K_\nu$ being the $K$-Bessel function of order $\nu$. The
$\varrho_j(n)$ are called the Fourier coefficients of $\psi_j$. In addition, we may
suppose that $\psi_j$ are simultaneous eigenfunctions of all Hecke operators with
corresponding eigenvalues $\tau_j(n)\in \B{R}$; that is, for each positive integer
$n$
$$
{1\over\sqrt{n}}\sum_{ad=n}\,\sum_{b\bmod d}\psi_j\left((az+b)/d\right)=
\tau_j(n)\psi_j(z),\quad z\in\B{H}.\eqno(1.5)
$$
We have, for any $m,\, n>0$,
$$
\tau_j(m)\tau_j(n)=\sum_{d|(m,n)}\tau_j\left(mn/d^2\right).\eqno(1.6)
$$
We may assume further that
$$
\psi_j\left(-\overline{z}\right)=\epsilon_j\psi_j(z),\quad
\epsilon_j=\pm1.\eqno(1.7)
$$
Then, the Hecke $L$-function associated with $\psi_j$ is 
defined by
$$
H_j(s)=\sum_{n=1}^\infty \tau_j(n)n^{-s},\quad \Re(s)>1.\
\eqno(1.8)
$$
This continues to an entire function, satisfying the functional equation
$$
H_j(s)=\chi_j(s)H_j(1-s),\eqno(1.9)
$$
with
$$
\eqalignno{
\chi_j(s)&=\epsilon_j\pi^{2s-1}{\Gamma\left({1\over2}\left(1-s+i\kappa_j
+{1\over2}(1-\epsilon_j)\right)\right)
\Gamma\left({1\over2}\left(1-s-i\kappa_j
+{1\over2}(1-\epsilon_j)\right)\right)
\over \Gamma\left({1\over2}\left(s+i\kappa_j
+{1\over2}(1-\epsilon_j)\right)\right)\Gamma\left({1\over2}\left(s-i\kappa_j
+{1\over2}(1-\epsilon_j)\right)\right)}\qquad&(1.10)\cr
&=2^{2s-1}\pi^{2(s-1)}\Gamma(1-s+i\kappa_j)
\Gamma(1-s-i\kappa_j)\left\{\epsilon_j\cosh(\pi\kappa_j)
-\cos(\pi s)\right\}.&(1.11)
}
$$
Since we have $\tau_j(n)\ll n^{1/4+\varepsilon}$ uniformly in $\psi_j$
(see [23, (3.1.18)]), the equation
$(1.9)$ implies that $H_j(s)$ is of polynomial growth 
with respect to both $s$ and $\kappa_j$ in any fixed vertical strip of the
$s$-plane.
\par
We shall need also holomorphic cusp forms over $\varGamma$, and
corresponding Hecke $L$-functions. Thus, if $\psi$ is holomorphic over $\B{H}$, 
vanishing at $i\infty$, and
$\psi(z)(dz)^k$ with a positive integer $k$ is $\varGamma$-invariant, then
we call it a holomorphic cusp form of weight $2k$. The space composed of all such
functions is a finite dimensional Hilbert space. We denote the dimension by
$\vartheta(k)$, and let $\left\{\psi_{j,k}:\,1\le j\le\vartheta(k)\right\}$ be a
corresponding orthonormal basis. Note that $\vartheta(k)=0$ for $k\le 5$. The
Fourier coefficient
$\varrho_{j,k}(n)$ of $\psi_{j,k}$ is defined by the expansion
$$
\psi_{j,k}(z)=\sum_{n=1}^\infty n^{k-1/2}\varrho_{j,k}(n)
\exp(2\pi iz),\quad z\in\B{H}.\eqno(1.12)
$$
We may assume that $\psi_{j,k}$ are simultaneous eigenfunctions of all Hecke
operators, so that there exist real numbers $\tau_{j,k}(n)$ such that
$$
{1\over\sqrt{n}}\sum_{ad=n}(a/d)^k\sum_{b\bmod d}
\psi_{j,k}\left((az+b)/d\right)=\tau_{j,k}(n)\psi_{j,k}(z),\quad z\in\B{H}.
\eqno(1.13)
$$
Then the Hecke $L$-function associated with $\psi_{j,k}$ is defined by
$$
H_{j,k}(s)=\sum_{n=1}^\infty \tau_{j,k}(n)n^{-s},\quad \Re(s)>1.\eqno(1.14)
$$
This continues to an entire function; and it satisfies the functional
equation
$$
H_{j,k}(s)=-
2^{2s-1}\pi^{2(s-1)}\Gamma\left(\txt{1\over2}-s+k\right)
\Gamma\left(\txt{3\over2}-s-k\right)
\cos(\pi s)H_{j,k}(1-s).\eqno(1.15)
$$
\medskip
Now, returning to our original subject, let $H$ be a particular function
among those $H_j$ and $H_{j,k}$. Comparing $(1.11)$ and $(1.15)$ with the
functional equation
$$
\zeta^2(s)=2^{2s-1}\pi^{2(s-1)}\Gamma^2\left(1-s\right)
\left(1-\cos(\pi s)\right)\zeta^2(1-s),\eqno(1.16)
$$
and invoking what is stated above about the size of
$\zeta^2\left({1\over2}+it\right)$, we might put $B(H,t)=t$; and an expected
subconvexity bound would be
$$
H\left(\txt{1\over2}+it\right)\ll t^{1/3+\varepsilon},\quad t\ge1.\eqno(1.17)
$$
In the case of holomorphic cusp forms, this was proved by A. Good [4] as a
corollary of an asymptotic formula for the mean square of 
$H\left({1\over2}+it\right)$, which he achieved by an appeal to the
spectral theory of real analytic automorphic functions (see also [22]).
An alternative and conceptually simpler proof, based solely on functional
properties of $H$ and its twists with additive characters (see $(8.8)$ below), 
was devised by the first named author [8]. His argument turned out to be
applicable also to the real analytic case, as shown by T. Meurman [19], yielding a
proof of $(1.17)$. Good's mean value result itself was later extended to this case
by the first named author [10], which implies $(1.17)$ in yet another way.
\medskip
With these developments in background, it should be desirable to have bounds
uniform in $\psi_j$. More precisely, $(1.9)$--$(1.11)$ suggest that we may choose
$B(H_j,t)=\kappa_j+t$ for $t\ge0$, and hence a hypothetical uniform subconvexity
bound would be $(1.1)$.
As a support, the first named author [13] showed recently that
$$
H_j\left(\txt{1\over2}+it\right)\ll t^{1/3+\varepsilon},\quad
t\gg\kappa_j^{3/2-\varepsilon},\eqno(1.18)
$$
which supersedes Meurman's
estimate, with respect to uniformity. This is in fact a consequence of the
following result on the spectral mean square (loc.cit):
$$
\sum_{K\le\kappa_j\le K+G}\alpha_j\left|H_j\left(\txt{1\over2}+it\right)\right|^2
\ll\left( GK+t^{2/3}\right)^{1+\varepsilon},\quad t\ge 0,
\; 1\le G\le K,\eqno(1.19)
$$
where $\alpha_j=|\varrho_j(1)|^2/\cosh(\pi\kappa_j)$.
Hence, when $t$ is relatively large, the bound $(1.1)$ holds indeed, in
view of the lower bound $\alpha_j\gg \kappa_j^{-\varepsilon}$ due to H. Iwaniec
[7]. The assertion $(1.19)$ has an essential relevance to our discussion in
Section 8, where a brief description of its proof is given.
\par
The real interest is, however, in the range 
$$
0\le t\le\kappa_j^{3/2},\eqno(1.20)
$$
since here the discrete quantity $\kappa_j$ seems to overwhelm the 
influence of the continuous parameter $t$. In this circumstance, what  A.
Ivi\'c [5] had achieved prior to $(1.18)$ was a breakthrough: He succeeded
in proving $(1.1)$ for
$t=0$ by a method quite different from those previously applied. His starting
point was an identity due to the second named author [23, Lemma 3.8] for the
spectral average
$$
\sum_{j=1}^\infty \alpha_j\tau_j(f)H_j^2\left(\txt{1\over2}\right)h(\kappa_j),
\eqno(1.21)
$$
where $h$ is a weight
function satisfying certain regularity and decay condition. As is
precisely presented in Lemma 3 below, this identity transforms the sum $(1.21)$
into a purely arithmetic expression involving, in particular, the divisor
function $d(n)$, which Ivi\'c could exploit effectively. His bound for
$H_j\left({1\over2}\right)$ is a corollary of the following result thus
obtained:
$$
\sum_{K\le\kappa_j\le K+G}\alpha_jH_j^3\left(\txt{1\over2}\right)
\ll GK^{1+\varepsilon}, \quad 1\le G\le K,\eqno(1.22)
$$
and the assertion $H_j\left({1\over2}\right)\ge0$ due to S. Katok and P. Sarnak
[15]. It should be remarked that a spectral sum of cubic powers of
$H_j\left({1\over2}\right)$ appeared for the first time in an explicit spectral
expansion of the weighted fourth moment of $\zeta\left({1\over2}+it\right)$ due to
the second named author [20] (see also [23, Chapter 4]).
Motivated by this advance with the cubic moment, the first named author [12]
turned to the fourth moment, establishing
$$
\sum_{K\le\kappa_j\le K+G}\alpha_jH_j^4\left(\txt{1\over2}\right)
\ll GK^{1+\varepsilon},\quad K^{1/3}\le G\le K.\eqno(1.23)
$$
The same identity for the sum $(1.21)$ played again a crucial r\^ole in his proof.
Also, as a new basic ingredient, a use was made of an explicit spectral
decomposition of the binary additive divisor sum
$$
D(f;W)=\sum_{n=1}^\infty d(n)d(n+f)W(n/f),\quad f\ge1,\eqno(1.24)
$$
due to the second named author [21] (see Lemma 5 below). It
should be stressed that $(1.23)$ proves Ivi\'c's bound for
$H_j\left({1\over2}\right)$ without the non-negativity assertion quoted after
$(1.22)$.
\medskip
Having stated this, it is now natural to investigate
the spectral fourth moment
$$
{\cal S}(G, K)=\sum_{K\le \kappa_j\le K+G}\alpha_j
\left|H_j\left(\txt{1\over2}+it\right)\right|^4,
\quad t\ge0,\; 1\le G\le K,\eqno(1.25)
$$
trying to retain the same bound as $(1.23)$ with uniformity in the parameter $t$.
Indeed, it gives rise to a proof of $(1.1)$:
\medskip
\noindent
{\bf Theorem 1.} {\it Let $K$ be sufficiently large, and
$$
G=(K+t)^{4/3}K^{-1+\varepsilon},\quad 0\le t\ll K^{3/2-\varepsilon}.\eqno(1.26)
$$
Then, we have
$$
{\cal S}(G,K)\ll GK^{1+\varepsilon}.\eqno(1.27)
$$
In particular, the bound $(1.1)$ holds uniformly for any $t\ge0$ and for
any real analytic cusp form $\psi_j$.
}
\medskip
\noindent
This embodies the main result of the present article. The second assertion follows
immediately from $(1.19)$ and $(1.27)$. For orientation, it should be remarked
that the estimate
$$
{\cal S}(G,K)\ll (K+t)^{2+\varepsilon},\eqno(1.28)
$$
with the same specification as in $(1.25)$, follows immediately from Lemmas 7 and
9 below.
\medskip
The proof of $(1.27)$ that we shall develop below is in principle an elaboration of
the argument in [12]. However, the prerequisite that the whole of our procedure
be  uniform in the parameter $t$ necessitates major changes of argument as well.
First of all, we are unable to exploit a peculiar property of the central
values of Hecke series, on which both [5] and [12] 
rely via the explicit formula for
$(1.21)$. Thus instead we appeal to the sum formula
of R.W. Bruggeman [1] and N.V. Kuznetsov [16] and in tandem to the sum formula of
Vorono{\"\i}. This is made at an earlier phase, i.e., Section 4, of the reduction
process, and causes already a considerable complication; nevertheless, it leads us
to an instance of the additive divisor sum $D(f;W)$. The subsequent procedure
is far more involved than the corresponding steps in [12], as will be
seen in Sections 5 and 6.  Moreover, only when $t$ is relatively small, i.e.,
$t\le K^{2/3}$, the end result thus reached is appropriate for an application of
the spectral large sieve (see Lemma 7 below) to produce what we desire. The analogy
with [12] ceases here.  For larger $t$ in the range $(1.26)$, the same
combination yields only an assertion short of $(1.27)$. 
Thence, we enter into the second phase of our discussion. That is about a spectral
hybrid mean value of Hecke series, an implement to extract $(1.27)$ out of the
aforementioned end result. This part might raise a particular interest, because a
significant contribution of holomorphic cusp forms takes place. It is thus
suggested that what we deal with in the present article is of quite a different
nature from any problem in analytic number theory to which  the spectral theory of
automorphic forms was applied,  e.g., the fourth moment of the Riemann
zeta-function, where the r\^ole of holomorphic cusp forms was in fact negligible.
\medskip
More precisely, the spectral hybrid mean value is concerned with the expression
$$
{\cal T}(K,t)=\sum_{K\le\kappa_j\le 2K}\alpha_jH_j^2\left(\txt{1\over2}\right)
\left|H_j\left(\txt{1\over2}+it\right)\right|^2.\eqno(1.29)
$$
\medskip
\noindent
{\bf Theorem 2.}\quad{\it 
We have, for any $K,\,t\ge0$,
$$
{\cal T}(K,t)\ll \left(K^2+t^{4/3}\right)^{1+\varepsilon}.\eqno(1.30)
$$
}
\par
\noindent
This is in fact an auxiliary result; thus it should be noted that no attempt is
made to prove the  best result obtainable by present day methods. The
proof developed in Section 7 starts with the explicit formula for $(1.21)$ and
follows to some extent the arguments of [12] and [23, Section 3.4]. We encounter an
additive divisor sum of the type
$$
D(f;\alpha,\beta;W)=\sum_{n=1}^\infty\sigma_\alpha(n)\sigma_\beta(n+f)W(n/f),
\quad f\ge1,\eqno(1.31)
$$
where $\sigma_\alpha(n)=\sum_{d|n}d^\alpha$, and in our situation $\alpha,\beta$
are complex. We can appeal to an explicit formula for this sum
due to the second named author [21]; however, the subsequent discussion is quite
subtle. We shall have two instances of $D(f;\alpha,\beta;W)$; and to deal
with the first, we require a counterpart of $(1.18)$--$(1.19)$ for holomorphic cusp
forms. This is precisely the peculiarity of our problem mentioned above. Thus,
uniformly for any
$\psi_{j,k}$
$$
H_{j,k}\left(\txt{1\over2}+it\right)\ll t^{1/3+\varepsilon},\quad
t\gg k^{3/2-\varepsilon}.\eqno(1.32)
$$
Also, under the same specification as in $(1.19)$,
$$
\sum_{K\le k\le K+G}\sum_{j=1}^{\vartheta(k)}\alpha_{j,k}
\left|H_{j,k}\left(\txt{1\over2}+it\right)\right|^2
\ll\left( GK+t^{2/3}\right)^{1+\varepsilon},\eqno(1.33)
$$
where $\alpha_{j,k}=8(4\pi)^{-2k-1}(2k-1)!|\varrho_{j,k}(1)|^2$.
The former is  of course a consequence of the latter together with 
an obvious analogue of the lower bound for $\alpha_j$.
On the other hand, with another instance of $D(f;\alpha,\beta;W)$, we require
instead $(1.19)$ in an analogous configuration. Therefore, the holomorphic and the
real analytic cusp forms stand at parity in our discussion of ${\cal T}(K,t)$.
\par
A proof of $(1.33)$ is given in the final section. It depends on an
observation about a crucial r\^ole played by the divisor function in our discussion
so far laid out. We are then led not only to $(1.33)$ but also to the
following counterpart of Theorem 1:
\medskip
\noindent
{\bf Theorem 3.}\quad{\it We have, under $(1.26)$,
$$
\sum_{K\le k\le K+G}\sum_{j=1}^{\vartheta(k)}\alpha_{j,k}\left|H_{j,k}
\left(\txt{1\over2}+it\right)\right|^4
\ll GK^{1+\varepsilon}.\eqno(1.34)
$$
In particular, the bound $(1.2)$ holds uniformly for any $t\ge0$ and for
any holomorphic cusp form $\psi_{j,k}$.
}
\medskip
With this, we look into the structure of our
argument, in order to envisage further extensions of our main result $(1.1)$; and
we come to a circle of problems on the size of
Rankin--Selberg $L$-functions. We shall indicate that our  method is capable of
yielding new results in such a generality as well.
\medskip
In passing, it should be added that the bounds $(1.1)$--$(1.2)$ could be stated
more uniformly, if we refer to basic terms from the theory of the
$\varGamma$-automorphic representations of the Lie group $\r{PSL}_2(\B{R})$,
which can be found in [3], for instance. Thus, we have
$$
H_V\left(\txt{1\over2}+it\right)\ll
\left(|\nu_V|+t\right)^{1/3+\varepsilon},
\eqno(1.35)
$$
uniformly for $t\ge0$ and for any Hecke invariant 
irreducible cuspidal representation $V$ with the spectral data
$\nu_V$, occurring in $L^2(\varGamma\backslash \r{PSL}_2(\B{R}))$. In the final
section we shall make a digression relevant to this aspect of our work.
\medskip
Throughout our discussion, the common symbol $\varepsilon$ plays a
basic r\^ole. Here we make precise our usage of them, in terms of a convention.
This is to avoid any confusion that might arise otherwise:
\medskip
\noindent
{\bf Convention 1.} The symbol $\varepsilon$ denotes
a sufficiently small positive parameter, which in general
takes different values at each occurrence. An
$\varepsilon_0>0$ could actually be fixed initially so
that a local value of $\varepsilon$ is an integral multiple of $\varepsilon_0$,
and  each inequality holds with an implied constant which depends solely on our
choice of $\varepsilon_0$. Thus, except being stated together with extra
dependencies, the notation $X\ll Y$, with $Y>0$, implies  that $|X|/Y$ is  bounded
by a constant depending on
$\epsilon_0$ at most, and $X\approx Y$ means that
$1\ll |X/Y|\ll 1$. It is implicit in our argument how to choose
multiples of $\varepsilon_0$ to have a particular inequality and a
specific reasoning valid. 
\medskip
Notations and conventions, including those in the above, are introduced where they
are needed for the first time, and will continue to be effective thereafter.
\vskip 1cm
\noindent
\centerline{\bf 2. Basic identities}
\bigskip
\noindent
Our proof of Theorem 1 is comprised of a series of various
transformations and approximations applied to spectral and arithmetic
objects. Here we collect identities which will give rise to fundamental
transformations of
${\cal S}(G,K)$ in later sections.
\medskip
\noindent
{\bf Lemma 1.} {\it Let $h(r)$ be even and regular in the strip 
$|\Im(r)|<{1\over4}+\varepsilon$, and there 
$|h(r)|\ll (1+|r|)^{-2-\varepsilon}$.
Put
$$
\eqalignno{
\widehat{h}(x)
&={2i\over\pi}\int_{-\infty}^\infty{rh(r)\over\cosh(\pi r)} J_{2ir}(x)dr\cr
&={2i\over\pi}\int_0^\infty{rh(r)\over\cosh(\pi
r)}\left(J_{2ir}(x)- J_{-2ir}(x)\right)dr,&(2.1)
}
$$
with $J_\nu$ being the $J$-Bessel function of order $\nu$. Then we have
$$
\eqalignno{
&\sum_{j=1}^\infty\alpha_j\tau_j(m)\tau_j(n)h(\kappa_j)=-{1\over\pi}
\int_{-\infty}^\infty{\sigma_{2ir}(m)\sigma_{2ir}(n)\over (mn)^{ir}
|\zeta(1+2ir)|^2}h(r)dr\cr
+&{\delta_{m,n}\over\pi^2}\int_{-\infty}^\infty r\tanh(\pi r)h(r)dr+
\sum_{\ell=1}^\infty{1\over\ell}S(m,n;\ell)\widehat{h}(4\pi\sqrt{mn}/\ell),
&(2.2)
}
$$
where $\delta_{m,n}$ is the Kronecker delta, and 
$$
S(m,n;\ell)=\sum_{\scr{q=1}\atop\scr{(q,\ell)=1}}^\ell\exp\left(2\pi i
(mq+n\tilde{q})/\ell\right),\quad q\tilde{q}\equiv1\bmod \ell,
\eqno(2.3)
$$
is a Kloosterman sum.
}
\medskip
\noindent
{\it Proof\/}. This is a refined version of the Spectral--Kloosterman sum formula
of Bruggeman and Kuznetsov. See [23, Section 2.6] for a proof.
\medskip
\noindent
{\bf Lemma 2.} {\it We have, for any integers $k,m,n\ge1$,
$$
\eqalignno{
\alpha_{j,k}\sum_{j=1}^{\vartheta(k)}&\tau_{j,k}(m)\tau_{j,k}(n)
={1\over2\pi^2}\delta_{m,n}(2k-1)\cr
&+{(-1)^k\over\pi}(2k-1)\sum_{\ell=1}^\infty
{1\over\ell}S(m,n;\ell)J_{2k-1}\left({4\pi}\sqrt{mn}/\ell\right).&(2.4)
}
$$
}
\medskip
\noindent
{\it Proof\/}. This is the sum formula of H. Petersson. A proof is given in
[23, Section 2.2].
\medskip
\noindent
{\bf Lemma 3.} {\it Let $h(r)$ be an even entire function satisfying 
$$
h\left(\pm\txt{1\over2}i\right)=0,\eqno(2.5)
$$
and
$$
h(r)\ll\exp(-\varepsilon|r|^2),\eqno(2.6)
$$
in any fixed horizontal strip. Put
$$
\eqalignno{
\Psi^+(x;h)&=2\pi\int_0^1\left\{y(1-y)(1+y/x)\right\}^{-{1/2}}\cr
&\times\int_{-\infty}^\infty rh(r)\tanh(\pi r)\left({y(1-y)\over
x+y}\right)^{ir}dr\,dy,&(2.7)
}
$$
and
$$
\eqalignno{
\Psi^-(x;h)=\int_0^\infty&\left\{\int_{(a)}x^s(y(y+1))^{s-1}
{\Gamma^2({1\over2}-s)\over \Gamma(1-2s)\cos(\pi s)}ds\right\}\cr
&\times\left\{\int_{-\infty}^\infty rh(r)
\left({y\over y+1}\right)^{ir}dr\right\}dy,
&(2.8)
}
$$
with $-{3\over2}<a<{1\over2}$, $a\ne-{1\over2}$, where $(a)$ 
is the vertical line
$\Re(s)=a$. Then we have, for any $f\ge1$,
$$
\sum_{j=1}^\infty\alpha_j\tau_j(f)H_j^2\left(\txt{1\over2}\right)
h(\kappa_j)=\sum_{\nu=1}^7{\e{H}}_\nu(f;h),
\eqno(2.9)
$$
where
$$
\eqalignno{
&{\e{H}}_1(f;h)=-i{2\over\pi^3}{d(f)\over\sqrt{f}}\cr
&\times \int_{-\infty}^\infty\left\{2\left(\gamma_E-\log(2\pi\sqrt{f})\right)
{\Gamma'\over\Gamma}\left(\txt{1\over2}+ir\right)
+\left({\Gamma'\over\Gamma}\left(\txt{1\over2}+ir\right)\right)^2
\right\}rh(r)dr,&(2.10)
}
$$
with the Euler constant $\gamma_E$, and
$$
{\e{H}}_2(f;h)={1\over\pi^3}\sum_{m=1}^\infty{m}^{-{1/2}}
d(m)d(m+f)\Psi^+\left({m/f};h\right),\eqno(2.11)
$$
$$
{\e{H}}_3(f;h)={1\over\pi^3}\sum_{m=1}^\infty(m+f)^{-{1/2}}
d(m)d(m+f)\Psi^-\left(1+{m/f};h\right),\eqno(2.12)
$$
$$
{\e{H}}_4(f;h)={1\over\pi^3}\sum_{m=1}^{f-1}m^{-{1/2}}
d(m)d(f-m)\Psi^-\left({m/f};h\right),\eqno(2.13)
$$
$$
{\e{H}}_5(f;h)=-{1\over2\pi^3}{d(f)\over\sqrt{f}}\Psi^-(1;h),\eqno(2.14)
$$
$$
{\e{H}}_6(f;h)=-{12\over\pi^2}f^{1/2}\sigma_{-1}(f)
h'\left(-\txt{1\over2}i\right),\eqno(2.15)
$$
$$
{\e{H}}_7(f;h)=-{1\over\pi}\int_{-\infty}^\infty f^{-ir}\sigma_{2ir}(f)
{{|\zeta({1\over2}+ir)|^4}\over
{|\zeta(1+2ir)|^2}}h(r)dr.\eqno(2.16)
$$
}
\medskip
\noindent
{\it Proof\/}. This is [23, Lemma 3.8]. 
Note that we have invoked the formulas [23, $(3.3.41)$ and $(3.3.45)$].
The decay condition on $h$ could be far less stringent than $(2.6)$.
\medskip
\noindent
{\bf Lemma 4.}\quad{\it Let $D(f;\alpha,\beta;W)$ be defined by $(1.31)$, where
$W$ is a smooth function supported compactly in the positive
reals, and $|\Re(\alpha)|,\,|\Re(\beta)|<\varepsilon$. Then we have
$$
D(f;\alpha,\beta;W)=\left\{D_r+D_d+D_h+D_c\right\}(f;\alpha,\beta;W),\eqno(2.17)
$$
where
$$
D_r(f;\alpha,\beta;W)=\int_0^\infty W(x)Y_f(x;\alpha,\beta)dx,\eqno(2.18)
$$
$$
\eqalignno{
\quad
&D_d(f;\alpha,\beta;W)=2(2\pi)^{\beta-1}f^{(\alpha+\beta+1)/2}
\sum_{j=1}^\infty\alpha_j\tau_j(f)\cr
&\times H_j\left(\txt{1\over2}(1-\alpha-\beta)\right)
H_j\left(\txt{1\over2}(1+\alpha-\beta)\right)
\left(\Psi_++\epsilon_j\Psi_-\right)
(i\kappa_j;\alpha,\beta;W), &(2.19)
}
$$
$$
\eqalignno{
\quad
&D_h(f;\alpha,\beta;W)=2(2\pi)^{\beta-1} f^{(\alpha+\beta+1)/2}
\sum_{k=1}^\infty\sum_{j=1}^{\vartheta(k)}(-1)^k
\alpha_{j,k}\tau_{j,k}(f)\cr
&\times H_{j,k}\left(\txt{1\over2}(1-\alpha-\beta)\right)
H_{j,k}\left(\txt{1\over2}(1+\alpha-\beta)\right)
\Psi_+\left(k-\txt{1\over2};\alpha,
\beta;W\right),\quad&(2.20)
}
$$
$$
\eqalignno{
&D_c(f;\alpha,\beta;W)=4(2\pi)^{\beta-2}f^{(\alpha+\beta+1)/2}\cr
&\times\int_{-\infty}^\infty f^{-i\kappa}\sigma_{2i\kappa}(f){
Z(i\kappa;\alpha,\beta)\over
\zeta(1+2i\kappa)\zeta(1-2i\kappa)}
\left(\Psi_++\Psi_-\right)(i\kappa;\alpha,\beta;W)
d\kappa. &(2.21)
}
$$
Here
$$
\eqalignno{
Y_f(x;\alpha,\beta)
&=\sigma_{1+\alpha+\beta}(f){\zeta(1+\alpha)\zeta(1+\beta)
\over\zeta(2+\alpha+\beta)}x^\alpha(x+1)^\beta\cr
&+f^\alpha\sigma_{1-\alpha+\beta}(f){\zeta(1-\alpha)\zeta(1+\beta)
\over\zeta(2-\alpha+\beta)}(x+1)^\beta\cr
&+f^\beta\sigma_{1+\alpha-\beta}(f){\zeta(1+\alpha)\zeta(1-\beta)
\over\zeta(2+\alpha-\beta)}x^\alpha\cr
&+f^{\alpha+\beta}\sigma_{1-\alpha-\beta}(f){\zeta(1-\alpha)\zeta(1-\beta)
\over\zeta(2-\alpha-\beta)},&(2.22)
}
$$
$$
\eqalignno{
Z(\xi;\alpha,\beta)&=\zeta\left(\txt{1\over2}(1-\alpha-\beta)+\xi\right)
\zeta\left(\txt{1\over2}(1+\alpha-\beta)+\xi\right)\cr
&\times\zeta\left(\txt{1\over2}(1-\alpha-\beta)-\xi\right)
\zeta\left(\txt{1\over2}(1+\alpha-\beta)-\xi\right),&(2.23)
}
$$
and
$$
\eqalignno{
&\qquad\Psi_+(\xi;\alpha,\beta;W)={1\over4\pi
i}\cos\left(\txt{1\over2}\pi\alpha\right)\int_{-i\infty}^{i\infty}
\cos(\pi s)\Gamma(s+\xi)
\Gamma(s-\xi)\cr
&\times\Gamma\left(\txt{1\over2}(1-\alpha-\beta)-s\right)
\Gamma\left(\txt{1\over2}(1+\alpha-\beta)-s\right)W^*
\left(s+\txt{1\over2}(\alpha+\beta+1)\right)ds,&(2.24)
 }
$$
$$
\eqalignno{
&\qquad\Psi_-(\xi;\alpha,\beta;W)={1\over4\pi
i}\cos\left(\pi\xi\right)\int_{-i\infty}^{i\infty}\sin\left(\pi\left(
s+\txt{1\over2}\beta\right)\right)\Gamma(s+\xi)
\Gamma(s-\xi)\cr
&\times\Gamma\left(\txt{1\over2}(1-\alpha-\beta)-s\right)
\Gamma\left(\txt{1\over2}(1+\alpha-\beta)-s\right)
W^*\left(s+\txt{1\over2}(\alpha+\beta+1)\right)ds,
&(2.25)
 }
$$
where $W^*$ is the Mellin transform of $W$, and the last integrals are such
that the path separates the poles of the first three factors in the integrand
and those of the remaining factors to the left and the right of the path,
respectively. 
}
\medskip
\noindent
{\it Proof\/}. This is asserted in [21, $(3.45)$--$(3.47)$], save for a
minor modification applied to the $D_h$ term. Also the formulas [21, $(3.42)$ 
and $(3.49)$] are invoked. The above condition on the real parts of $\alpha,\,
\beta$ is imposed only for the sake of convenience, and thus by no means
essential. In fact, the explicit formula $(2.17)$ holds for all complex
$\alpha,\,\beta$ in the context of analytic continuation. In [21] it is
implicitly assumed that $W$ is real valued, but in fact the argument there
allows us to drop it; hence in the above $W$ can be complex valued.
\medskip
\noindent
{\bf Lemma 5.}\quad{\it Let $D(f;W)$ be defined by $(1.24)$, with $W$
being as in the previous lemma, and 
let
$$
\eqalignno{
Y_f(u)&=(\log u)\log(u+1)+\left(c-\log f+2{\sigma_1'\over\sigma_1}(f)\right)
\log(u(u+1))\cr
&+(c-\log
f)^2-4\left({\zeta'\over\zeta}\right)'\!(2)
+4{\sigma_1'\over\sigma_1}(f)(c-\log
f) +{\sigma_1''\over\sigma_1}(f),&(2.26) }
$$
where $\sigma_\xi^{(\nu)}=(d/d\xi)^\nu\sigma_\xi$, and
$c=2\gamma_E-2(\zeta'/\zeta)(2)$.
Then we have
$$
D(f;W)=\left\{D_r+D_d+D_h+D_c\right\}(f;W),\eqno(2.27)
$$
where
$$
\eqalignno{
D_r(f;W)&={6\over\pi^2}\sigma_1(f)\int_0^\infty Y_f(u)W(u)du,&(2.28)\cr
D_d(f;W)&=f^{1/2}\sum_{j=1}^\infty\alpha_j\tau_j(f)H_j^2
\left(\txt{1\over2}\right)\Psi(i\kappa_j;W),&(2.29)\cr
D_h(f;W)&=f^{1/2}\sum_{\scr{k=6}\atop\scr{2|k}}^\infty
\sum_{j=1}^{\vartheta(k)}\alpha_{j,k}\tau_{j,k}(f)H_{j,k}^2
\left(\txt{1\over2}\right)
\Psi\left(k-\txt{1\over2};W\right),&(2.30)\cr
D_c(f;W)&={f^{1/2}\over\pi}
\int_{-\infty}^\infty f^{-i\kappa}\sigma_{2i\kappa}(f)
{|\zeta({1\over2}+i\kappa)|^4\over|\zeta(1+2i\kappa)|^2}
\Psi(i\kappa;W)d\kappa.&(2.31)
}
$$
Here
$$
\eqalignno{
\Psi(\xi;W)={1\over2}&\int_0^\infty
\Re\Bigg\{\left(1-{1\over\sin(\pi\xi)}\right)
{\Gamma^2({1\over2}+\xi)\over\Gamma(1+2\xi)}\cr
&\times 
{}_2F_1\left(\txt{1\over2}+\xi,\txt{1\over2}+\xi;
1+2\xi;-{1/u}\right)u^{-{1/2}-\xi}\Bigg\}W(u)du,&(2.32)
}
$$
with the Gaussian hypergeometric function ${}_2F_1$.
}
\medskip
\noindent
{\it Proof\/}. This is a corollary of the last lemma. See [21] for details. Note
that here it is invoked that $H_j\left({1\over2}\right)=0$ if $\epsilon_j=-1$, 
and $H_{j,k}\left({1\over2}\right)=0$ if $k$ is odd, as the functional equations
$(1.9)$ and $(1.15)$ imply, respectively.
\medskip
Lemmas 3--5 should be compared with corresponding assertions claimed by Kuzne\-tsov
[17]. We add also that in the light of [3] the explicit formula $(2.17)$ could be
derived directly from the spectral structure of $L^2(\varGamma\backslash
\r{PSL}_2(\B{R}))$, that is, without appealing to the spectral theory of
sums of Kloosterman sums on which [17] and [21] rely.
\vskip 1cm
\centerline{\bf 3. Basic inequalities}
\bigskip
\noindent
In the present section we shall prepare those implements which are crucial in our
approximation procedures pertaining to estimations of our key objects. Asymptotics in this context will be supplied
mostly by the saddle point method. The proof of Lemma 7 below furnishes typical
instances which could be referred to at later applications of the method. 
Note that Convention 1 is always in force hereafter.
\par
To facilitate the relevant reasoning and in fact the whole of our discussion,  the
following formulation of the treatment of off-saddle integrals will turn out to be
highly instrumental:
\medskip
\noindent
{\bf Lemma 6.} {\it Let $A$ be a smooth function compactly supported in a finite
interval $[a,b]$; and assume that there exist two quantities $A_0$, $A_1$ such that
for each integer $\nu\ge0$ and for any $x$ in the interval
$$
A^{(\nu)}(x)\ll A_0A_1^{-\nu}.\eqno(3.1)
$$ 
Also, let $B$ be a function which is real-valued on $[a,b]$, and regular 
throughout the complex domain composed of all points within the distance 
$\rho$ from the interval; and assume that there exists a quantity $B_1$ such that
$$
0<B_1\ll |B'(x)|\eqno(3.2)
$$
for any point $x$ in the domain. Then we have, for each fixed integer $P\ge0$,
$$
\int_{-\infty}^\infty A(x)\exp\left(iB(x)\right)dx\ll
A_0(A_1B_1)^{-P}(1+A_1/\rho)^P(b-a).
\eqno(3.3)
$$ 
\/}
\par
\noindent
{\it Proof\/}. With a multiple application of integration by parts, 
we see that the integral is equal to
$$
i^P\int_a^b\left[\left({\cal
D}_B\right)^P\!A\right](x)\exp\left(iB(x)\right)dx,\eqno(3.4)
$$
where ${\cal D}_B$ is the operator $g\mapsto (g/B')'$. We have
$$
\eqalignno{
\left({\cal D}_B\right)^P\!A(x)&=\sum_{\nu_1+\cdots+\nu_P\le P}a(\nu_1,\ldots,
\nu_P)A^{(P-\nu_1-\cdots-\nu_P)}(x)\cr
&\times\left({1\over B'(x)}\right)^{(\nu_1)}\cdots
\left({1\over B'(x)}\right)^{(\nu_P)},&(3.5)
}
$$
with certain constants $a(\nu_1,\ldots,\nu_P)$. The assumption $(3.2)$
gives $\left(1/B'(x)\right)^{(\nu)}\ll B_1^{-1}\rho^{-\nu}$ on $[a, b]$ via
Cauchy's integral formula for derivatives. Thus,  $(3.1)$ implies that in $(3.4)$
$$
\left({\cal D}_B\right)^P\!A(x)\ll A_0(A_1B_1)^{-P}\sum_{\nu_1+\cdots+\nu_P\le P}
(A_1/\rho)^{\nu_1+\cdots+\nu_P},\eqno(3.6)
$$
from which $(3.3)$ follows.
\medskip
\noindent
{\bf Lemma 7.} {\it Let $1\le G\le K$ and $N\ge1$. Then we have, for any
complex vector $\{a(n)\}$,
$$
\sum_{K\le\kappa_j\le K+G}\alpha_j\left|\sum_{N\le n\le 2N}\tau_j(n)
a(n)\right|^2\ll(GK+N)(KN)^\varepsilon\sum_{N\le n\le 2N}|a(n)|^2.
\eqno(3.7)
$$
}
\medskip
\noindent
{\it Proof\/}. This version of the spectral large sieve of Iwaniec [6] is due to
the first named author [11] (see [23, Theorem 3.3] for a refinement). Here we
shall show a new approach to $(3.7)$.  A truncation procedure in our argument,
i.e., $(3.12)$ below, will turn out to be  fundamental for our discussion of
${\cal S}(G,K)$ that starts in the next section. It should be noted that 
smooth and compactly supported weights attached to integers could be avoided in the
present proof proper; their use is made rather for the sake of later purpose.
\par 
Obviously we may assume that
$K^\varepsilon\ll G\ll K^{1-\varepsilon}$, with the basic parameter $K$ that is
larger than a constant depending solely on $\varepsilon_0$. 
The case $N\gg K^{1/\varepsilon}$ can be settled by an application of
a duality principle and the theory of Rankin zeta-functions (see [23, pp.\
137--138]). Thus, we  may assume also that $N\ll K^{1/\varepsilon}$.  With this,
let
$$
h(r)=K^{-2}\left(r^2+\txt{1\over4}\right)\left[\exp\left(-((r-K)/G)^2\right)+
\exp\left(-((r+K)/G)^2\right)\right].\eqno(3.8)
$$
It suffices to prove that 
$$
\sum_{j=1}^\infty\alpha_j\left|\sum_{N\le n\le 2N}\phi(n)
\tau_j(n)a(n)\right|^2h\left(\kappa_j\right)\eqno(3.9)
$$
is bounded by the right side of $(3.7)$, where $\phi$ is an arbitrary real-valued
smooth function which is supported in $\left[N,2N\right]$ and 
$\phi^{(\nu)}(y)\ll N^{-\nu}$ for each $\nu\ge0$. Expand out the square, take the
spectral sum inside, and apply $(2.2)$. The contribution to
$(3.9)$ of the first term on the right of $(2.2)$ is negative, and can be
discarded. That of the second term is obviously absorbed into the right side of
$(3.7)$.  Then, let
${\cal A}$ be the part of $(3.9)$ corresponding to the sum of Kloosterman sums on
the right of
$(2.2)$.  In ${\cal A}$, the sum over $\ell$ can be truncated to $1\le\ell\ll
K^{1/\varepsilon}$, under Convention 1. This can be seen by shifting the contour
of the first integral in $(2.1)$ to $\Im(r)=-1$. In fact, we have 
$$
\widehat{h}(x)\ll\int_{-\infty}^\infty {(|r|+1)h(r)\over\cosh(\pi r)}
|J_{2+2ir}(x)|dr\ll Gx^2/K,\eqno(3.10)
$$
via Poisson's formula
$$
J_\nu(x)={(x/2)^\nu\over\sqrt{\pi}\Gamma(\nu+{1\over2})}
\int_{-1}^1(1-y^2)^{\nu-{1/2}}\cos(xy)dy,\eqno(3.11)
$$
which is valid for $x>0$, $\Re(\nu)>-{1\over2}$. 
\par
We shall show that the remaining part of ${\cal A}$ could be truncated to
$$
1\le\ell\ll N(GK)^{-1}\log K.\eqno(3.12)
$$
To this end, we invoke the representation
$$
J_{2ir}(x)-J_{-2ir}(x)={2\over\pi i}\sinh(\pi r)\Re
\int_{-\infty}^\infty
\exp\left(ix\cosh u-2iru\right)du\eqno(3.13)
$$
(see the formula (12) on [26, p.\ 180]), which we use with $x=4\pi\sqrt{mn}/\ell$.
When $|u|>\log^2K$, we perform integration by parts with respect to the factor
$\exp(ix\cosh u)$, getting, for $r\ge0$,
$$
\eqalignno{
J_{2ir}(x)-J_{-2ir}(x)={2\over\pi i}&\sinh(\pi r)\Re
\int_{-\log^2K}^{\log^2K}
\exp\left(ix\cosh u-2iru\right)du\cr
&+O\left((r+1)\exp\left(\pi r-\txt{1\over2}\log^2K\right)\right).
&(3.14)
}
$$
Thus, via the second expression in $(2.1)$, we obtain, after a rearrangement,
$$
\eqalignno{
\widehat{h}(x)
&=K^{-2}\Re\int_{-\infty}^{\infty}R(u)\eta(u)\exp(-(Gu)^2)
\exp\left(ix\cosh u-2iKu\right)du\cr
&\hskip 2cm+O\left(\exp\left(-\txt{1\over3}\log^2K\right)\right),
&(3.15)
}
$$
where $R$ is a certain polynomial on $u,G,K$, and 
$\eta(u)$ is a smooth weight such that
$\eta(u)=1$ for $|u|\le (\log K)/G$, and $=0$ for $|u|\ge 2(\log K)/G$ as well as
$\eta^{(\nu)}(u)\ll (G/\log K)^\nu$ for each $\nu\ge0$. In fact, the
expression follows first with the range $|u|\le\log^2K$ but without the weight;
then the truncation to $|u|\le (\log K)/G$ can be imposed; and the result is
modified as $(3.15)$. With this, we assume temporarily that
$x\ll GK/\log K$. Then, Lemma 6 is applicable to the last integral, with
$A_0=GK^3$, $A_1=G^{-1}\log K$, $B_1=K$, $\rho\approx G^{-1}$. We find that
$\widehat{h}(x)$ is negligibly small. Hence, we may restrict ourselves to $x
\gg GK/\log K$, which is the same as $(3.12)$. In passing, we note that
$$
R(u)\eta(u)=4\pi^{-3/2}GK^3\left(1+O(K^{-\varepsilon})\right)\eta(u).\eqno(3.16)
$$
\par
We may thus equip $\cal A$ with weights $\phi_0(\ell)$, where $\phi_0$ is a smooth
function such  that $\phi_0(y)=0$ for $y\le0$,
$=1$ for $1\le y\le L_0$, $=0$ for $y\ge 2L_0$;  here $L_0$ is a dyadic number 
such that $L_0\ll (N\log K)/(GK)$ and the error  thus caused is negligible. We
replace $\phi_0(y)$ by a sum of smooth $\phi_1(y;L)$ such that it is
supported in $[L/2,2L]$ with a dyadic $L\le L_0$, and $\phi_1^{(\nu)}(y)\ll
L^{-\nu}$ for each $\nu\ge0$. Hence, it suffices to deal
with
$$
\sum_{m=1}^\infty\sum_{n=1}^\infty\phi(m)\phi(n)a(m)\overline{a(n)}
\sum_{\ell=1}^\infty{\phi_1(\ell)\over\ell}S(m,n;\ell)\widehat{h}
\left(4\pi\sqrt{mn}/\ell\right),\eqno(3.17)
$$
where $\phi_1(y)=\phi_1(y;L)$.
Obviously we require
$$
GKL\ll N\log K.\eqno(3.18)
$$
\par
Then, we shall show that provided $(3.18)$ we have
$$
\eqalignno{
\widehat{h}(x)&={4\over\pi}\sqrt{2\over x}GKR_1(x)\exp\left(-\left({Gu_0}
\right)^2\right)\cr
&\times\cos\left(x\cosh u_0-2Ku_0+\txt{1\over4}\pi\right)
+O\left(K^{-1/\varepsilon}\right),&(3.19)
}
$$
where $\sinh u_0=2K/x$ or $u_0=\log\left(2K/x+\sqrt{1+4(K/x)^2}\right)$; 
and
$$
R_1(x)=\sum_{\nu}b_\nu(G,K)x^\nu,\eqno(3.20)
$$
with a finite number of terms; $b_\nu(G,K)x^\nu\ll K^{-\varepsilon}$
for $\nu\ne0$ and $b_0(G,K)=1$. Note that any regular function of
$u_0$ is a power series of $K/x$, and could be approximated
by a polynomial on $K/x$ with an arbitrary accuracy; this is relevant to
our reasoning below.
\par
Before entering into the proof of $(3.19)$, let us make a useful observation: What
matters in estimating $(3.17)$ is to fix the leading term in
$(3.20)$; that is, in the asymptotic evaluation of $\widehat{h}(x)$, which we are
about to develop using the saddle point method, we may restrict our attention
to the main term, provided it is clear
that the argument yields, in fact, an expansion of the type $(3.19)$ with
$(3.20)$. This  is due to the fact that the contribution to $(3.17)$ of the term 
$b_\nu(G,K)x^\nu$ ($\nu\ne0$) of $(3.20)$ could be dealt with in just the same way
as that of the term with $\nu=0$,  because it corresponds to a change of weight:
$$
\eqalignno{
&\phi(m)\phi(n)\phi_1(\ell)\cr
&\mapsto b_\nu(G,K)(M/L)^{\nu}
[\phi(m)(m/M)^{\nu/2}][\phi(n)(n/M)^{\nu/2}][\phi_1(\ell)(\ell/L)^{-\nu}].
&(3.21)
}
$$
The new weight thus obtained is of the same type as the original but only smaller,
as the factor $b_\nu(G,K)(M/L)^{\nu}$ is $\ll K^{-\varepsilon}$ for $\nu\ne0$. 
\par
With this, we consider $(3.15)$; then the above observation
allows us to treat instead the simpler
$$
J(x)=4\pi^{-3/2}GK \int_{-\infty}^{\infty}\exp(-(Gu)^2)
\exp\left(ix\cosh u-2iKu\right)du,\eqno(3.22)
$$
since it will be clear from our discussion that the factor $R(u)$ gives rise
to a factor of the type $R_1(x)$; and since the weight $\eta(u)$ can obviously be
removed. We now apply the saddle point method to $(3.22)$, which is routine but
better to be performed with some details because of our later purpose. 
Thus, $u_0$ is the saddle point, and $(3.12)$ or $(3.18)$ gives $u_0\ll K/x\ll
G^{-1}\log K$.  We put $u=v+\xi\exp\left({1\over4}\pi i\right)$.
We move the last contour to $C_-+C_0+C_+$, where 
$C_-=\{u: v<u_0, \xi=-\xi_0\}$, $C_0=\{u: v=u_0, -\xi_0\le\xi\le\xi_0\}$,
$C_-=\{u: u_0<v, \xi=\xi_0\}$, with an obvious orientation and
$\xi_0= x^{-1/(2+\varepsilon)}$. Accordingly, the dissection
$J(x)=\left\{J^{(-)}+J^{(0)}+J^{(+)}\right\}(x)$ follows. Note that 
we have $\xi_0\ll (GK^{\varepsilon})^{-1}$,
because of $(3.18)$, $G\ll K^{1-\varepsilon}$, and Convention 1. In particular, 
$\exp(-(Gu)^2)\ll\exp\left(-{1\over2}(Gv)^2\right)$ 
on the new contour. We have also
$$
\eqalignno{
x\cosh u&-2Ku=x\cosh v-2Kv+\xi\exp\left(\txt{1\over4}\pi i\right)
\left(x\sinh v-2K\right)\cr
&+{i\over2}x\xi^2\cosh v +x\sum_{j=3}^\infty{1\over j!}\left(\xi
\exp\left(-\txt{1\over4}\pi i\right)\right)^j
\cosh\left(v+\txt{1\over2}j\pi i\right).&(3.23)
}
$$
This implies that 
$\Im(x\cosh u-2Ku)>{1\over3}x\xi_0^2\cosh v$ on $C_\pm$, since  $\pm(x\sinh
v-2K)>0$  throughout $C_\pm$. Hence we have
$(J^{(-)}+J^{(+)})(x)\ll K\exp\left(-{1\over3}x\xi_0^2\right)$, which is
negligible. On the other hand, we have
$$
\eqalignno{
J^{(0)}(x)&=4\pi^{-3/2}GK \exp\left(ix\cosh u_0-2iKu_0+\txt{1\over4}\pi
i\right)\cr &\times\int_{-\xi_0}^{\xi_0}\exp\left(-(Gu)^2\right)
\exp\left(-\txt{1\over2}x\xi^2\cosh u_0+i\Sigma\right)d\xi,&(3.24)
}
$$
where $\Sigma$ is the last term of $(3.23)$ with $v=u_0$. Since $\Sigma\ll
x\xi_0^3$, the factor $\exp(i\Sigma)$ can be replaced by a polynomial on
$\Sigma$ with a negligible error; and the power series in $\Sigma$ is
truncated with the same effect. Also, $\exp\left(-(Gu)^2\right)$ is replaced
by $\exp\left(-(Gu_0)^2\right)$ times a polynomial in a similar fashion.
This and applications of integration by parts give that $J^{(0)}(x)$ is equal
to a multiple by a factor of the type $R_1(x)$ of
$$
\eqalignno{
4\pi^{-3/2}&GK \exp\left(-(Gu_0)^2\right)
\exp\left(ix\cosh u_0-2iKu_0+\txt{1\over4}\pi
i\right)\cr &\times\int_{-\xi_0}^{\xi_0}
\exp\left(-\txt{1\over2}x\xi^2\cosh u_0\right)d\xi
+O\left(K^{-1/\varepsilon}\right),&(3.25)
}
$$
which leads us to the assertion $(3.19)$. 
\par
Now we return to the estimation of $(3.17)$. As observed above, it suffices
to consider the contribution of $(3.19)$ with $R_1(x)$ being replaced by $1$. Then,
we put
$$
\theta_{*}(s)=\int_0^\infty \theta(x)\exp\left(-(Gu_0)^2
\right)\exp\left(ix\cosh
u_0-2iKu_0+\txt{1\over4}\pi i\right)x^{s-1} dx,\eqno(3.26)
$$
where $\theta(x)$ is a smooth function which is equal to $1$ over $[2\pi
N/L,8\pi N/L]$, supported in $[\pi N/L,10\pi N/L]$, and 
$\theta^{(\nu)}(x)\ll x^{-\nu}$ 
for each $\nu\ge0$. We find that the estimation of the part of $\cal A$ under
consideration is, via the Mellin inversion of $(3.26)$ and the definition $(2.3)$,
reduced to that of
$$
GK L^{-{1/2}}\sum_{L\le\ell\le2L}
\sum_{\scr{q=1}\atop\scr{(q,\ell)=1}}^\ell\int_{-\infty}^\infty
\left|\sum_{N\le n\le 2N} {\phi(n)a(n)\over n^{(1/2+iv)/2}}\exp\left(2\pi
i{qn/\ell}\right)
\right|^2|\theta_{*}(iv)|dv.\eqno(3.27) 
$$
This is, by the hybrid large sieve inequality (see [23, Lemma 3.11]),
$$
\ll GK(LN)^{-{1/2}}\sum_{V\ge 1}
(N+L^2(V+1))\sup_{V\le|v|\le2V}|
\theta_{*}(iv)|\sum_{N\le n\le2N}|a(n)|^2,\eqno(3.28)
$$
where $V$ runs over dyadic numbers.
\par
To bound $\theta_{*}(iv)$, we note that
$$
{d\over dx}\left(v\log x+x\cosh u_0-2Ku_0\right)=v/x+
\sqrt{1+4(K/x)^2}.\eqno(3.29)
$$
The saddle point $x_0$ of the integral $(3.26)$ with $s=iv$ is close to $-v$, and
has to be inside the support of $\theta$, since otherwise $\theta_*(iv)$
could be seen to be negligibly small by Lemma 6, with $A_0=L/N$, $A_1=N/L$,
$B_1=1+VL/N$, $\rho\approx N/L$. In particular, we may assume that
$$
V\approx N/L.\eqno(3.30)
$$ 
With this, we shall further estimate
$\theta_{*}(iv)$ by the saddle point method. 
Thus, let $\rho_0=(N/L)^{2/5}$, and
let $\theta_1$ is a smooth function such that 
$\theta_1(x)=1$ for $|x-x_0|\le
x_0/\rho_0$, and $=0$ for $|x-x_0|\ge 2x_0/\rho_0$; as well as
$\theta_1^{(\nu)}(x)\ll (x_0/\rho_0)^{-\nu}$ for each $\nu\ge0$. Also, let
$\theta_2=\theta-\theta_1$. Then, Lemma 6 implies that 
$(\theta_2)_*(iv)$ is negligibly small; in fact, this results
with the specification $A_0=L/N$, $A_1=(N/L)^{3/5}$, $B_1=(N/L)^{1/5}$, 
$\rho\approx(N/L)^{3/5}$. Thus, via the Taylor expansion of the 
integrand of $(3.26)$  at $x=x_0$, we have that 
$$
\eqalignno{
\theta_*(iv)=&\exp\left(-(Gu_0)^2\right)
\exp\left(ix_0\cosh u_0-2iKu_0+\txt{1\over4}\pi i\right)x_0^{iv-1}\cr
&\times\int_{-\infty}^\infty R_2(x)
\theta_1(x)\exp\left(\txt{1\over2}iX(x-x_0)^2\right)dx
+O(K^{-1/\varepsilon}),&(3.31)
}
$$
where $u_0$ is specialised with $x=x_0$; $X\approx L/N$ is the derivative of
$(3.29)$ at $x=x_0$; and $R_2(x)$ is analogous to $(3.20)$. This integral is
divided into two parts according as $|x-x_0|\le X^{-1/2}$ and otherwise. The first
part is bounded trivially, and the second after an application of integration by
parts. We obtain 
$$
\theta_*(iv)\ll (N/L)^{-{1/2}}.\eqno(3.32)
$$
Being inserted into the part of $(3.28)$ corresponding to $(3.30)$, 
this gives rise to the assertion $(3.7)$. 
\medskip
\noindent
{\bf Lemma 8.} {\it 
With the same specifications as in the previous lemma, we have
$$
\sum_{K\le k\le K+G}\sum_{j=1}^{\vartheta(k)}
\alpha_{j,k}\left|\sum_{N\le n\le 2N}
\tau_{j,k}(n)a(n)\right|^2\ll(GK+N)(KN)^\varepsilon\sum_{N\le n\le 2N}|a(n)|^2.
\eqno(3.33)
$$
}
\medskip
\noindent
{\it Proof\/}. To show this counterpart of
$(3.7)$, we put  $h_1(r)=\exp(-\left((K-r)/G\right)^2)$, multiply the inner sum on
the left of $(3.33)$ by the factor $h_1(k)$, and sum over all integers $k\ge1$. 
Then a use of $(2.4)$ leads us to a sum of Kloosterman sums as in the previous
proof, but with $\widehat{h}(x)$ being replaced by
$$
\sum_{k=1}^\infty (-1)^k(2k-1)h_1(k)J_{2k-1}(x).\eqno(3.34)
$$
By the Poisson integral
$$
J_{2k-1}(x)={(-1)^k\over\pi}\int_{-{1\over2}\pi}^{{1\over2}\pi}
\sin\left((2k-1)u-x\cos u)\right)du\eqno(3.35)
$$
and the Poisson sum formula, one may see that $(3.34)$ can be replaced by
$$
{2\over\sqrt{\pi}}GK\int_{-(\log K)/G}^{(\log K)/G}
\sin\left(2Ku-x\cos u\right)\exp\left(-(Gu)^2\right)du.\eqno(3.36)
$$
With this, the rest of the proof is analogous to the above.
\bigskip
\noindent
{\bf Lemma 9.} {\it Let $K$ be a large positive parameter. Let 
$$
K^\varepsilon\ll G\ll K^{1-\varepsilon},\quad 0\le t\ll
K^{1/\varepsilon},\eqno(3.37)
$$
and put
$$
T={1\over4\pi^2}(K+t)\left(|K-t|+G\right).\eqno(3.38)
$$
Then, uniformly for all cusp forms $\psi_j$ with
$$
|K-\kappa_j|\ll G,\eqno(3.39)
$$
we have
$$
\eqalignno{
H_j^2&\left(\txt{1\over2}+it\right)
\ll (\log K)^2\sum_{M\le 2T}\sum_{q=1}^\infty{1\over q}\cr
&\times\int_{\gamma^{-1}-i\gamma^2}^{\gamma^{-1}+i\gamma^2}
\left|\sum_{m=1}^\infty
\phi(q^2m;M)d(m)\tau_j(m)m^{-\xi-{1/2}-it}\right|
|d\xi|,&(3.40)
}
$$
with $\gamma=\log^2K$ and dyadic numbers $M$.
Here the smooth function $\phi(y;M)$ depends solely on the interval 
$\left[M/2,2M\right]$, in which it is supported and $\phi^{(\nu)}(y;M)
\ll M^{-\nu}$, with the implied constant depending only on $\nu$.
\/}
\medskip
\noindent
{\it Proof.\/} We consider the integral
$$
{\cal{R}}={1\over{2\pi{i}\gamma}}\int_{(3)}H_j^2\left(\xi+\txt{1\over2}
+it\right)T^\xi\Gamma(\xi/\gamma)d\xi.\eqno(3.41)
$$
Since $(1.6)$  gives
$$
H_j^2(s)=\zeta(2s)\sum_{n=1}^\infty d(n)\tau_j(n)n^{-s}, \quad \Re(s)>1,
\eqno(3.42)
$$
we have 
$$
{\cal{R}}=\sum_{n\le\alpha T}\check{\tau}_j(n)n^{-{1/2}-it}\exp(-(n/T)^\gamma)
+O(e^{-K}).\eqno(3.43)
$$
where $2\le\alpha\le4$ is arbitrary, and
$$
\check{\tau}_j(n)=\sum_{q^2|n}d(n/q^2)\tau_j(n/q^2).\eqno(3.44)
$$
Shifting the path of $(3.41)$ to $\left(-{1\over2}\gamma\right)$ and
recalling the functional equation $(1.10)$, we get
$$
{\cal{R}}=H_j^2\left(\txt{1\over2}+it\right)+\sum_{n=1}^\infty
\check{\tau}_j(n)n^{-{1/2}+it}{\cal{R}}_j(n),\eqno(3.45)
$$
where
$$
\eqalignno{
{\cal{R}}_j(n)&={\pi^{4it}\over2\pi i\gamma}\int_{(-{1\over2}\gamma)}
\left({\Gamma\left({1\over2}\left({1\over2}-\xi-it+i\kappa_j
+{1\over2}(1-\epsilon_j)\right)\right)
\over \Gamma\left({1\over2}\left({1\over2}+\xi+it-i\kappa_j+
{1\over2}(1-\epsilon_j)\right)\right)}\right)^2\cr
&\times\left({\Gamma\left({1\over2}\left({1\over2}-\xi-it-i\kappa_j
+{1\over2}(1-\epsilon_j)\right)\right)\over
\Gamma\left({1\over2}\left({1\over2}+\xi+it+i\kappa_j+
{1\over2}(1-\epsilon_j)\right)\right)}\right)^2
(\pi^4 nT)^\xi\Gamma\left({\xi/\gamma}\right)d\xi.&(3.46)
}
$$
By Stirling's formula this integrand is
$$
\ll(16\pi^4nT)^{-\gamma/2}\big(|\xi+i(t-\kappa_j)|
|\xi-i(\kappa_j+t)|\big)^\gamma\exp(-\pi|\xi|/(2\gamma));\eqno(3.47)
$$
and thus
$$
{\cal{R}}_j(n)\ll (n/T)^{-\gamma/2}.\eqno(3.48)
$$
In fact, when $|\xi|<\gamma^2$ we see that
the factor $\big(|\xi+i(\kappa_j-t)| |\xi-i(\kappa_j+t)|\big)^{\gamma}$ is
$\ll (4\pi^2T)^\gamma$ in view of $(3.39)$, and when
$|\xi|\ge\gamma^2$ the integrand itself is negligible due to the factor
$\exp(-\pi|\xi|/(2\gamma))$. The estimate $(3.48)$ allows us to truncate
the sum in $(3.45)$ to $n\le \alpha T$. 
In this way, we have, uniformly for all $\psi_j$
satisfying $(3.39)$,
$$
\eqalignno{
H_j^2\left(\txt{1\over2}+it\right)=
&\sum_{n\le{\alpha T}}\check{\tau}_j(n)n^{-{1/2}-it}\exp(-(n/T)^\gamma)\cr
-&\sum_{n\le{\alpha T}}\check{\tau}_j(n)n^{-{1/2}+it}{\cal R}_j(n)
+O\left(K^{-1}\right). &(3.49)
}
$$
We equip both sums with smooth and compactly supported
weights in much the same way as performed preceding $(3.17)$; here 
the parameter $\alpha$ plays a r\^ole. Then the first sum in
$(3.49)$ is readily seen to be bounded by the right side of
$(3.40)$.  As to the second sum, we modify
${\cal{R}}_j$ by shifting the path in $(3.46)$ to $(-\gamma^{-1})$, and take
the sum over $n$ inside the integral. Considering the absolute value of the
resulting integrand, we may eliminate the $\Gamma$-factors of $(3.46)$ except for
$\Gamma(\xi/\gamma)$. This gives rise to $(3.40)$.
\medskip
\noindent
{\bf Corollary to Lemma 9.}\quad{\it Let 
$K^\varepsilon\ll G\ll K^{1-\varepsilon}$. Then we have 
$$
\sum_{K\le\kappa_j\le K+G}\alpha_j
\left|H_j\left(\txt{1\over2}+it\right)\right|^4\ll
GK^{1+\varepsilon},\eqno(3.50)
$$
uniformly for $|K-t|\ll G$. 
}
\par
\medskip
\noindent
{\it Proof.\/} This follows immediately from a combination of Lemmas 7 and 9.
\medskip
With these preparations we shall start our discussion
of ${\cal S}(G,K)$, in the next section. Technically it is a layered
application of those approximation--estimation procedures employed in the proof
of Lemma 7. To avoid excessive repetitions of details, we introduce
\medskip
\noindent
{\bf Convention 2.} All subsequent approximations are to hold with the
basic parameter $K$ that is assumed to be
larger than a quantity depending solely on
$\varepsilon_0$.  With this, let ${\cal X}$ be a particular object that we need to
bound. Suppose that an expression $\cal Y$ comes up in a relevant discussion, and
we have an approximation ${\cal Y}={\cal Y}_0+{\cal Y}_1+O({\cal Z})$, in which 
${\cal Y}_0$ is dominant, ${\cal Y}_1$ oscillates in the same mode as
${\cal Y}_0$, while $\cal Z$ contributes negligibly to ${\cal X}$. Then, the
notation ${\cal Y}\sim {\cal Y}_0$ indicates an actuation of a procedure in which
the treatment of ${\cal Y}_1$ is a repetition of that of
${\cal Y}_0$ and the replacement of
${\cal Y}$ by ${\cal Y}_0$ causes no differences in bounding ${\cal X}$.
\medskip
\noindent
For instance, in the proof of Lemma 7, the polynomial factor $R_1(x)$ of $(3.19)$
is essentially irrelevant to the estimation of $(3.17)$; and this could be denoted
as $R_1(x)\sim 1$. More drastically, as we shall do in the sequel, this economy of
reasoning could have been applied to $(3.17)$ from the very beginning of the
proof, as $(3.21)$ endorses. We shall employ devices analogous to $(3.21)$,
without mentioning persistently.
\vskip 1cm
\centerline{\bf 4. Reduction}
\bigskip
\noindent
We begin our discussion of ${\cal S}(G,K)$.
We assume that $K$ is as in Convention 2, and that $(1.26)$
holds. Note that $G\ll K^{1-\varepsilon}$ under Convention 1.
\medskip
Let $h$ be defined by $(3.8)$ but with the present specification of the
parameters.  Then, by Lemma 9, it suffices to treat
$$
\sum_{j=1}^\infty\alpha_j\left|\sum_{m=1}^\infty
\phi_0(m){d(m)\tau_j(m) m^{-1/2-it}}\right|^2h(\kappa_j)\eqno(4.1)
$$
where $\phi_0(x)=\phi(q^2x;M)x^{-\xi}$ with $\xi$ and $\phi(q^2x;M)$ as
in $(3.40)$, while $T$ is defined by $(3.38)$ with the present $G$. 
Thus, $\phi_0(x)$ is smooth, compactly supported accordingly, 
and $\phi_0^{(\nu)}(x)\ll \left((\log K)^4/x\right)^\nu$. 
\par
We proceed just in the same way as in the proof of Lemma 7. 
What is essential for our purpose is to bound the Kloosterman-sum part of $(4.1)$
thus obtained. In view of $(3.12)$, we may assume that the corresponding
truncation has already been performed to the present sum over the moduli of
Kloosterman sums.  Thus, more specifically, we shall consider
$$
\eqalignno{
{\cal S}_1
=\sum_{m=1}^\infty&\sum_{n=1}^\infty\phi(m)\phi(n)d(m)d(n)
(mn)^{-{1/2}}(m/n)^{it}\cr
&\times\sum_{\ell=1}^\infty{\phi_1(\ell)\over\ell}
S(m,n;\ell)\widehat{h}\left(4\pi\sqrt{mn}/\ell\right),
&(4.2)
}
$$
which corresponds to $(3.17)$.
Here $\phi$ and $\phi_1$ are real-valued smooth functions, which are 
compactly supported in $[M/2, 2M]$ and $[L/2, 2L]$, respectively, with
$$
GKL/\log K\ll M\ll T.\eqno(4.3)
$$ 
Also, we have $\phi^{(\nu)}(x)\ll\left((\log K)^4/M\right)^{\nu}$ and
$\phi_1^{(\nu)}(x)\ll L^{-\nu}$ for each $\nu\ge0$. Note that the present $M$
stands for $M/q^2$ in $(4.1)$. The symbols $M$, $L$, $\phi$, and $\phi_1$ will
retain the current specifications till the end of Section 6. 
We have
$$
{\cal S}(G,K)\ll GK^{1+\varepsilon}+|{\cal S}_1|K^\varepsilon.\eqno(4.4)
$$
In the sequel, we shall modify or transform the sum ${\cal S}_1$ in several steps.
The most significant contribution will be denoted by ${\cal S}_\nu$, $\nu=2,3,4$;
accordingly the estimation of ${\cal S}_1$ is reduced to the same for
${\cal S}_4$.
\medskip
We return to the second expression in $(2.1)$. We note that the integration
can be restricted to 
$$
|r-K|\ll G\log K,\eqno(4.5)
$$ 
because of the uniform bound $|J_{2ir}(x)|\le
\left(\cosh(2\pi r)\right)^{1/2}$ which follows from $(3.11)$. 
We then evaluate the integral $(3.13)$ asymptotically;
we require $x=4\pi\sqrt{mn}/\ell$ to appear in $(4.2)$, i.e.,
$\phi(m)\phi(n)\phi_1(\ell)\ne0$. Obviously, we may proceed in
much the same way as $(3.22)$--$(3.25)$, and get
$$
\eqalignno{
J_{2ir}(x)-J_{-2ir}(x)&\sim 
{\sqrt{2}\over i\sqrt{\pi x\cosh u_1}}e^{\pi r}
\cos\left(x\cosh u_1-2ru_1+\txt{1\over4}\pi \right)\cr
&\sim {1\over i}{\sqrt{2\over\pi x}}e^{\pi r}
\cos\left(\omega(r,x)+\txt{1\over4}\pi \right),&(4.6)
}
$$
where $x\sinh u_1=2r$, and 
$$
\omega(r,x)=x(1-2(r/x)^2).\eqno(4.7)
$$
That is, the left side of $(4.6)$ is asymptotically equal, within a negligible
error, to the right side multiplied by a factor similar to $R_1(x)$ defined at
$(3.20)$.  Here we have used the facts that
$x\cosh u_1=x+2r^2/x+O(r^4/x^3)$,
$ru_1=2r^2/x+O(r^4/x^3)$, and $r^4/x^3\ll K^4/(GK/\log K)^3\ll
K^{-\varepsilon}$ because of $(1.26)$ and $(4.3)$. 
\medskip
Hence, by Convention 2, it suffices to consider the expression
$$
\eqalignno{
{2^{3/2}\over\pi^2}
&\int_{K-G\log K}^{K+G\log K}rh(r)\Bigg\{\sum_{m=1}^\infty
\sum_{n=1}^\infty\phi(m)\phi(n)d(m)d(n)
(mn)^{-{3/4}}(m/n)^{it}\cr
&\times\sum_{\ell=1}^\infty{\phi_1(\ell)\over\sqrt{\ell}}S(m,n;\ell)
\cos\left(\omega(r,4\pi\sqrt{mn}/\ell)+\txt{1\over4}\pi\right)\Bigg\}dr.
&(4.8)
}
$$
This reduces the estimation of ${\cal S}_1$
to that of
$$
\eqalignno{
{\cal S}_2&=\sum_{m=1}^\infty{\phi(m)d(m)\over m^{{3/4}-it}}
\sum_{\ell=1}^\infty{\phi_1(\ell)
\over\sqrt\ell}\sum_{\scr{q=1}\atop\scr{(q,\ell})=1}^\ell 
\exp(2\pi iqm/\ell)\cr
&\times\sum_{n=1}^\infty{\phi(n)d(n)\over n^{{3/4}+it}}
\exp(2\pi i\tilde{q}n/\ell)
\exp\left(\delta_1i\omega(r,4\pi\sqrt{mn}/\ell)\right),&(4.9)
}
$$
where $\delta_1=\pm1$, $q\tilde{q}\equiv1\bmod\ell$. Note that ${\cal S}_2$ is a
function of $r$. We have, via $(4.4)$,
$$
{\cal S}(G,K)\ll GK^{1+\varepsilon}\left(1+\sup_r|{\cal S}_2|\right),\eqno(4.10)
$$
where $r$ is in the range $(4.5)$.
\medskip
To transform ${\cal S}_2$, we apply the sum formula of Vorono\"{\i} 
(see, e.g.,\ [8, Theorem 1.7]) to the inner-most sum of $(4.9)$: Thus, it is equal
to
$$
\eqalignno{
{1\over\ell}\int_0^\infty&\left(\log y+2\gamma_E-2\log\ell\right)p(y)dy\cr
+&{1\over\ell}\sum_{n=1}^\infty d(n)\int_0^\infty\Big\{
4\exp\left(2\pi inq/\ell\right)K_0\left(4\pi\sqrt{ny}/\ell\right)\cr
&\qquad\qquad-2\pi\exp\left(-2\pi
inq/\ell\right)Y_0\left(4\pi\sqrt{ny}/\ell\right)\Big\}p(y)dy, &(4.11)
}
$$
where $K_0$ and $Y_0$ are Bessel functions in the
notation of [26], and 
$$
p(y)=\phi(y)y^{-{3/4}-it}
\exp\left(\delta_1i\omega(r,4\pi\sqrt{my}/\ell)\right).\eqno(4.12)
$$
For the sake of a later purpose, we stress that $(4.11)$ is a simple
consequence of the functional equation for the Hecke-Estermann zeta-function
(see, e.g.,\ [23, Lemma 3.7]):
Let
$$
D(s,\xi; q/\ell)=\sum_{n=1}^\infty\sigma_\xi(n)
\exp\left(2\pi inq/\ell\right)
n^{-s},\quad(q,\ell)=1.
\eqno(4.13)
$$
Then,
$$
\eqalignno{
&D(s,\xi;q/\ell)=2(2\pi)^{2s-\xi-2}\ell^{\xi-2s+1}\Gamma(1-s)\Gamma(1+\xi-s)\cr
&\times\left\{\cos\left(\txt{1\over2}\pi\xi\right)D(1-s,-\xi;\tilde{q}/\ell)
-\cos\left(\pi\left(s-\txt{1\over2}\xi\right)\right)
D(1-s,-\xi;-\tilde{q}/\ell)\right\},
&(4.14)
}
$$
which is actually equivalent to the automorphy of the real-analytic Eisenstein
series of weight $0$ over $\varGamma$.
\medskip
The leading term of $(4.11)$ is negligible by Lemma 6. In fact, we may 
set $\rho\approx M$; and in the relevant domain of $y$
$$
{d\over dy}\left(-t\log y
+\delta_1\omega(r,4\pi\sqrt{my}/\ell)\right)=-{1\over y}
\left(t+2\pi{\delta_1\over\ell}
\sqrt{my}+\delta_1{r^2\ell\over4\pi\sqrt{my}}\right).
\eqno(4.15)
$$
Here we have $r^2\ell/\sqrt{my}\ll K^2L/M\ll G^{-1}K\log K$ by $(4.3)$ and
$(4.5)$; and $\sqrt{m|y|}/\ell\gg GK/\log K\gg tK^\varepsilon$ by $(1.26)$.
Thus, Lemma 6 works with $A_0=(\log K)/M^{3/4}$, $A_1=M/(\log K)^4$, 
$B_1=GK/(M\log K)$, $\rho\approx M$; note that we have used the bound 
$\phi^{(\nu)}(y)\ll ((\log K)^4/M)^\nu$. This confirms our
claim. Also, the part of $(4.11)$ which contains the Bessel function
$K_0$ is negligible, because of the exponential decay of the function. 
As to the $Y_0$-part, we use the fact that
$Y_0(x)\sim (2/(\pi x))^{1/2}\sin\left(x-{1\over4}\pi\right)$ according to the
formula $(4)$ on [26, p.\ 199].  Thus, the main part of $(4.11)$ is 
$$
-\sqrt{2\over\ell}\sum_{n=1}^\infty{d(n)\over
n^{1/4}}\exp(-2\pi inq/\ell)\int_0^\infty y^{-{1/4}}p(y)\sin\left(
4\pi\sqrt{ny}/\ell-\txt{1\over4}\pi\right)dy.\eqno(4.16)
$$
Inserting this into $(4.9)$, we see that instead of ${\cal S}_2$ we may
deal with
$$
{\cal S}_3=\sum_{m=1}^\infty{\phi(m)d(m)\over m^{{3/4}-it}}
\sum_{\ell=1}^\infty{\phi_1(\ell)\over\ell}\sum_{n=1}^\infty{d(n)\over
n^{1/4}} c_\ell(m-n)I(\ell,m,n;\delta_1,\delta_2),\eqno(4.17)
$$
where $c_\ell$ is the Ramanujan sum $\bmod\,\ell$, and
$$
I(\ell,m,n;\delta_1,\delta_2)=\int_0^\infty\phi(y)y^{-1-it}\exp\left(\delta_1
i\omega(r,4\pi\sqrt{my}/\ell)+4\pi\delta_2i\sqrt{ny}/\ell\right)dy,\eqno(4.18)
$$
with $\delta_2=\pm1$. By Convention 2, we have, in place of $(4.10)$,
$$
{\cal S}(G,K)\ll GK^{1+\varepsilon}
\left(1+\sup_r|{\cal S}_3|\right).\eqno(4.19)
$$
\medskip
We apply Lemma 6 to the last integral. If
$\delta_1=\delta_2$, then the integral is similar to the leading term of $(4.11)$,
and can be discarded. Thus, hereafter we shall have $\delta_2=-\delta_1$. We may
set $\rho\approx M$ again, and in the relevant domain of $y$ we have
$$
\eqalignno{
&{d\over dy}\left(-t\log y
+\delta_1\omega(r,4\pi\sqrt{my}/\ell)
-4\pi\delta_1\sqrt{ny}/\ell\right)\cr
&={1\over y}\left(-t+2\pi\delta_1
{\sqrt{y}\over\ell}(\sqrt{m}-\sqrt{n})+\delta_1{r^2\ell\over4\pi\sqrt{my}}
\right).&(4.20)
}
$$
Let us assume that $|m-n|\gg L\left(t+{K^2L/M}\right)K^\varepsilon$. 
Then, throughout the domain we have
$(\sqrt{|y|}/\ell)\left|\sqrt{m}-\sqrt{n}\right|\gg
\left(t+{r^2\ell/\sqrt{m|y|}}\right)K^\varepsilon$. Hence,
Lemma 6 works with $A_0=1/M$, $A_1=M/(\log K)^4$, 
$B_1=\left(t+{K^2L/M}\right)K^\varepsilon/M$, $\rho\approx M$. Note
that $A_1B_1\gg K^{\varepsilon}$;
in fact, if $t\ge1$ then this is obvious, and otherwise $(3.38)$ and $(4.3)$ yield
the same.  Thus, $(4.18)$ is negligibly small, provided the above lower bound
for $|m-n|$. In other words, we may proceed with the truncation
$$
m-n\ll L\left(t+{K^2L\over M}\right)K^\varepsilon
\ll {M\over GK}\left(t+{K\over G}\right)K^\varepsilon.\eqno(4.21)
$$
\par
Let us settle the case $m=n$; that is, we are dealing with the
diagonal part of ${\cal S}_3$:
$$
\sum_{n=1}^\infty{\phi(n)d^2(n)\over n^{1-it}}\sum_{\ell=1}^\infty
{\phi_1(\ell)\over\ell}\varphi(\ell)I(\ell,n,n;\delta_1,-\delta_1),
\eqno(4.22)
$$
where $\varphi$ is Euler's totient function. We have
$$
I(\ell,n,n;\delta_1,-\delta_1)=\int_{1}^\infty
\phi(y)y^{-1-it}\exp\left(-\delta_1i{r^2\ell\over2\pi\sqrt{ny}}\right)dy,
\eqno(4.23)
$$
since $\phi(y)=0$ for $y\le 1$.
We can assume that $L\gg
K^\varepsilon$, for otherwise $(4.22)$ could obviously be ignored. 
Then, consider the situation $t\ll K^\varepsilon$; in particular
$T\approx K^2$, and $M\ll K^2$ by $(4.3)$. We may apply Lemma 6 to $(4.23)$,
with $A_0=1/M$, $A_1=M$, $B_1=K^2L/M$, $\rho\approx M$, 
since $K^2L/M\gg L\gg tK^\varepsilon$
under Convention 1. That is, this case can be ignored.
Let us move to the situation $t\gg K^\varepsilon$. We shall employ an argument
based on Mellin inversion; one may use the saddle point method as well. With the
Mellin transform $\phi^*$ of $\phi$, $(4.23)$ is equal to
$$
\eqalignno{
&{1\over2\pi i}
\int_{(\varepsilon)}\phi^*(s)\int_{1}^\infty
y^{-1-it-s}\exp\left(-\delta_1i{r^2\ell\over2\pi\sqrt{ny}}\right) dy\,ds\cr
=&{1\over2\pi i}
\int_{(\varepsilon)}\phi^*(s)\left\{\int_0^\infty-\int_0^{1}\right\}
y^{-1-it-s}\exp\left(-\delta_1i{r^2\ell\over2\pi\sqrt{ny}}\right) 
dy\,ds.\qquad
&(4.24)
}
$$
Note that $\phi^*(s)$ is of fast decay with respect to $s$ in any fixed vertical
strip. The double integral arising from
the last finite integral vanishes, as it can be seen by performing integration by
parts in the $y$-integral and exchanging the order of integration.
We have 
$$
I(\ell,n,n;\delta_1,-\delta_1)=
{1\over\pi i}
\int_{(\varepsilon)}\phi^*(s)\Gamma(2(s+it))e^{-\delta_1\pi i(s+it)}
\left({r^2\ell\over2\pi\sqrt{n}}\right)^{-2(s+it)}ds,\eqno(4.25)
$$
which converges absolutely.
Thus, the inner sum of $(4.22)$ can be written as
$$
\eqalignno{
-{1\over2\pi^2}\int_{(2)}&\phi_1^*(s_1)\int_{(\varepsilon)}\left({r^4\over
4\pi^2n}\right)^{-s-it}{\zeta(s_1+2(s+it))\over\zeta(s_1+2(s+it)+1)}\cr
&\times\phi^*(s)\Gamma(2(s+it))e^{-\delta_1\pi i(s+it)}ds\,ds_1,&(4.26)
}
$$
where $\phi_1^*$ is the Mellin transforms of $\phi_1$. This double integral
can be truncated to $|s|,\,|s_1|\ll K^\varepsilon$. Moving the $s_1$-contour to
the vertical line $(\varepsilon)$, we do not encounter any poles under Convention
1, and find that $(4.26)$ is $\ll K^\varepsilon$, which settles the present case.
Hence the diagonal part of ${\cal S}_3$ can be ignored.
\medskip
We turn to the non-diagonal part of ${\cal S}_3$:
$$
\eqalignno{
&\sum_{f\le f_0}\sum_{\ell=1}^\infty {\phi_1(\ell)c_\ell(f)\over\ell}
\sum_{n=1}^\infty{\phi(n)d(n)d(n+f)\over n^{{3/4}-it}(n+f)^{1/4}}
I(\ell,n,n+f;\delta_1,-\delta_1)\cr
+&\sum_{f\le f_0}\sum_{\ell=1}^\infty {\phi_1(\ell)c_\ell(f)\over\ell}
\sum_{n=1}^\infty{\phi(n+f)d(n)d(n+f)\over (n+f)^{{3/4}-it}n^{1/4}}
I(\ell,n+f,n;\delta_1,-\delta_1)\cr
=&{\cal S}_3^-+{\cal S}_3^+,&(4.27)
}
$$ 
say, where by $(4.21)$
$$
f_0\ll L\left(t+{K^2L/ M}\right)K^\varepsilon.\eqno(4.28)
$$
We have got instances of the additive divisor sum.
Let us put
$$
\eqalignno{
W_-(u)&=\phi(fu)u^{-{3/4}+it}(u+1)^{-{1/4}}I_-(f,\ell,u),\cr
W_+(u)&=\phi(f(u+1))(u+1)^{-{3/4}+it}u^{-{1/4}}I_+(f,\ell,u),&(4.29)
}
$$
where
$$
I_\pm(f,\ell,u)=\int_0^\infty{\phi(y)\over y^{1+it}}\exp\left(
\pm\delta_1i{4\pi\sqrt{fy}\over\ell(\sqrt{u}+\sqrt{u+1})}-\delta_1i{r^2\ell\over
2\pi\sqrt{fy(u+a_\pm)}}\right)dy,\eqno(4.30)
$$
with $a_\pm={1\over2}(1\pm1)$. Then $(4.27)$ can be written as
$$
{\cal S}_3^{\pm}=\sum_{f\le f_0}{1\over
f^{1-it}}\sum_{\ell=1}^\infty {\phi_1(\ell)c_\ell(f)\over\ell}
\sum_{m=1}^\infty d(m)d(m+f)W_\pm(m/f).\eqno(4.31)
$$
\par
Let us consider ${\cal S}_3^-$. With the change of variable
$v=4\pi\ell^{-1}\left(u+{1\over2}\right)^{-{1/2}}\sqrt{fy}$, 
we rewrite it as
$$
2(4\pi)^{2it}\sum_{f\le f_0}{1\over
f^{1-2it}}\sum_{\ell=1}^\infty {\phi_1(\ell)c_\ell(f)\over \ell^{1+2it}}
\sum_{m=1}^\infty d(m)d(m+f)W_-^{(1)}(m/f),\eqno(4.32)
$$
where
$$
W_-^{(1)}(u)=\phi(fu)u^{-{3/4}}(u+1)^{-{1/4}}
\left(1+{1\over2u}\right)^{-it}I_-^{(1)}(f,\ell,u),\eqno(4.33)
$$
with
$$
\eqalignno{
I_-^{(1)}&(f,\ell,u)=\int_0^\infty\phi\left({(\ell
v)^2\over16\pi^2f}\left(u+\txt{1\over2}\right)\right)\cr
&\times\exp\left(-\delta_1i{v\sqrt{u+{1\over2}}\over
\sqrt{u}+\sqrt{u+1}}-2\delta_1 i{r^2\over v\sqrt{u\left(u+{1\over2}\right)}}
\right)v^{-1-2it}dv.&(4.34)
}
$$
Here we could introduce the truncation
$$
u\approx M/f,\quad v\approx f/L,\eqno(4.35)
$$
in which the former is obvious, and the latter is due to the presence of the
$\phi$-factor in $(4.34)$. 
\par
We are going to simplify $W_-^{(1)}$ under Conventions 1 and 2.
To this end we note first that
$$
u^{-{3/4}}(u+1)^{-{1/4}}\left(1+{1\over2u}\right)^{-it}\sim
{1\over u}\exp\left(-{it\over2u}+{it\over8u^2}\right),\eqno(4.36)
$$
since $t/u^3\ll t(f/M)^3\ll t(GK)^{-3}(t+K/G)^3K^\varepsilon\ll
(GK)^{-3}(t+K)^4K^\varepsilon$, which is $\ll K^{-\varepsilon}$ because of
$(1.26)$. Also
$$
{v\sqrt{u+{1\over2}}\over\sqrt{u}+\sqrt{u+1}}\sim{v\over2}\left(1+
{1\over32 u^2}\right),\eqno(4.37)
$$
since $v/u^3\ll L^3M^{-3}(t+K/G)^4\ll (GK)^{-3}(t+K/G)^4K^\varepsilon\ll
K^{-\varepsilon}$. Further,
$$
{r^2\over v\sqrt{u\left(u+{1\over2}\right)}}\sim{r^2\over
uv}\left(1-{1\over4u}\right),\eqno(4.38)
$$
since $r^2/(u^3v)\ll K^2(L/M)^3(t+K/G)^2\ll(GK)^{-3}(t+K)^4K^\varepsilon
\ll K^{-\varepsilon}$.
\par
This leads us to
$$
{\cal S}_4=\sum_{f=1}^\infty{\phi_2(f)\over f^{1-2it}}
\sum_{\ell=1}^\infty{\phi_1(\ell)c_\ell(f)\over \ell^{1+2it}}
\sum_{m=1}^\infty d(m)d(m+f)X(m/f).\eqno(4.39)
$$
Here $\phi_2$ is a smooth function compactly supported in $[F/2,2F]$, with
$$
F\ll L(t+K/G)K^\varepsilon,\eqno(4.40)
$$
as it follows from $(4.28)$; and
$$
X(u)={1\over u}\int_0^\infty\xi(f,\ell,u,v)
\exp\left(iY\right){dv\over v^{1+2it}},\eqno(4.41)
$$
with
$$
\xi(f,\ell,u,v)=\phi(fu)\phi\left({(\ell v)^2\over
16\pi^2f}u\right)\eqno(4.42)
$$
and
$$
Y=-{t\over 2u}\left(1-{1\over4u}\right)-{1\over2}\delta_1v\left(1+
{1\over32 u^2}\right)-2\delta_1{r^2\over
uv}\left(1-{1\over4u}\right).\eqno(4.43)
$$
The $(4.42)$ depends on the fact that $\phi\left({((\ell v)^2/
16\pi^2f})(u+{1\over2})\right)\sim\phi\left({((\ell v)^2/
16\pi^2f})u\right)$. We note also that $\phi_2^{(\nu)}(x)\ll F^{-\nu}$
for each $\nu\ge0$, as usual.
\medskip
The transformation of ${\cal S}_3^+$ is
analogous. In fact we end up with the same expression as ${\cal S}_4$ except
for  the change of the definition $(4.43)$ into
$$
{t\over 2u}\left(1-{3\over4u}\right)+{1\over2}\delta_1v\left(1+
{1\over32 u^2}\right)-2\delta_1{r^2\over
uv}\left(1-{3\over4u}\right).\eqno(4.44)
$$
This should imply
that the discussion of ${\cal S}_3^+$ can be done with unessential alterations to
that of ${\cal S}_4$. Hence it suffices to treat ${\cal S}_4$; that is,
we have, in place of $(4.19)$, 
$$
{\cal S}(G,K)\ll GK^{1+\varepsilon}
\left(1+\sup_r|{\cal S}_4|\right),\eqno(4.45)
$$
with a minor abuse of reasoning. For a later convenience we note that
$(4.35)$ can be stated as
$$
u\approx M/F,\quad v\approx F/L,\eqno(4.46)
$$
with $(1.26)$, $(3.38)$, $(4.3)$, $(4.4)$, $(4.5)$, $(4.40)$ being provided.
The assertion $(4.45)$ is naturally dependent on a reasoning similar to $(3.21)$.
\vskip 1cm
\centerline{\bf 5. Lower range}
\bigskip
\noindent
We have reduced the estimation of ${\cal S}(G,K)$, a spectral object, to 
that of ${\cal S}_4$, an arithmetic object. With this, we now return to the
spectra. That is to say, we apply Lemma 5 to ${\cal S}_4$:
$$
{\cal S}_4=S_r+S_d+S_h+S_c,\eqno(5.1)
$$
in an obvious correspondence to the terms on the right of $(2.27)$. In the
present and the subsequent sections we shall deal with $S_r$ and $S_d$ in two
ranges of the parameter $t$. The parts $S_h$ and $S_c$ will be briefly treated;
they are analogous to $S_d$ and turn out to be negligible. 
\medskip
As vaguely indicated in Introduction, the range of $t$ is
divided into three sections according to the size of the spectral data under
consideration. This is rendered in the
division
$$
0\le t\le K^{2/3},\quad K^{2/3}\le t\le K^{3/2},\quad K^{3/2}\le t. \eqno(5.2)
$$
We call these intervals the lower, the intermediate, and the upper ranges,
respectively. The bound $(1.1)$ for the upper range follows from the spectral mean
square $(1.19)$; the difference caused by
those factors $\kappa_j^{\pm\varepsilon}$ and $K^{\pm\varepsilon}$ is immaterial
for our current discussion. Thus we consider the remaining two ranges. 
In the present section we shall deal with the lower range, or more precisely
we shall consider the situation
$$
0\le t\ll K^{1+\varepsilon}/G,\quad G \approx K^{1/3+\varepsilon}.\eqno(5.3)
$$
Note that consequentially we have $T\approx K^2$ and $M\ll K^2$.
As a matter of fact, this case has already been settled in the
announcement article [14], and thus could be skipped. Nevertheless, there is a
certain need to fill in some details missing in [14], and above all what we are
about to develop here should motivate effectively the reasoning in the next
section, where we shall treat the intermediate range. Also, we shall depart from
[14]  in a few technical aspects in order to show a variety of available methods.
\medskip
It should, however, be noted that our division of the range of $t$ is not 
imperative; a refinement of the argument of the next section should
make the present section redundant, at the cost of accessibility.
\medskip
Thus, we assume $(5.3)$, and consider the spectral expansion $(5.1)$.
Then we observe that
$$
Y\sim Q=-{1\over2}\delta_1v-2\delta_1{r^2\over uv}.\eqno(5.4)
$$
In fact, in $(4.43)$ we have $t/u\ll tF/M \ll t(K/G)(L/M)K^\varepsilon\ll
(t/G^2)K^\varepsilon$, as $(4.46)$, $(4.40)$, $(4.3)$
successively imply. Also, $v/u^2\ll F^3/(LM^2)\ll
(K/G)^3(L/M)^2K^\varepsilon\ll K^{1+\varepsilon}/G^5$ and $ r^2/(u^2v)\ll
K^2FL/M^2\ll K^{2+\varepsilon}(K/G)(L/M)^2\ll K^{1+\varepsilon}/G^3$.
Thus terms of $Y$, save for those two on the right of $(5.4)$, are all
negligible under Conventions 1 and 2.
\medskip
We shall consider ${\cal S}_4$ with $Y$ being replaced by $Q$. 
We begin with $S_r$. We may naturally take into account
only the leading term on the right of $(2.26)$, since the other terms are
treated similarly. On noting
$(\log u)\log (u+1)\sim(\log u)^2$, we need to bound the expression 
$$
\sum_{f=1}^\infty{\phi_2(f)\sigma_1(f)\over f^{1-2it}}
\sum_{\ell=1}^\infty{\phi_1(\ell)c_\ell(f)\over \ell^{1+2it}}
\int\!\int\xi(f,\ell,u,v)(\log u)^2\exp\left(iQ\right)
{dudv\over uv^{1+2it}},\eqno(5.5)
$$
where the range of integration is indicated by $(4.46)$. 
\par
Performing the change of variable $u\mapsto w/v$, we consider instead
$$
\sum_{f=1}^\infty{\phi_2(f)\sigma_1(f)\over f^{1-2it}}
\sum_{\ell=1}^\infty{\phi_1(\ell)c_\ell(f)\over \ell^{1+2it}}
\int\!\int\xi(f,\ell,w/v,v)(\log w/v)^2
\exp\left(iQ_1\right){dvdw\over v^{1+2it}w},\eqno(5.6)
$$
where
$$
v\approx F/L,\quad w\approx M/L,\eqno(5.7)
$$
and
$$
Q_1=-{1\over2}\delta_1v-2\delta_1{r^2\over w}.\eqno(5.8)
$$
We integrate with respect to $w$ first. Lemma 6 is applicable with $A_0=(\log
K)^2L/M$, $A_1=M/(L\log^4K)$, $B_1=(KL/M)^2$, $\rho\approx M/L $.  Thus the
$w$-integral of $(5.6)$ is $\ll (K^2L/(M\log^4K))^{-P}\log^2K$; and
$(5.6)$ can be discarded if $L\gg K^\varepsilon$, since we have $M\ll K^2$. 
That is, we may assume that $L\ll K^\varepsilon$. Then, obviously the case $F\ll
K^\varepsilon$ can be ignored. In particular, we may
assume that  $F/L\gg K^\varepsilon$ as well, under Convention 1. With this, we
integrate, in $(5.6)$, with respect to $v$ first. We have
$$
{\partial\over\partial v}(Q_1-2t\log v)=
-{1\over2}\delta_1-2{t\over v}.\eqno(5.9)
$$
In the case $\delta_1=1$, Lemma 6 can be applied with $A_0=(\log K)^2L/F$, 
$A_1=F/(L\log^4K)$, $B_1=1$, $\rho\approx F/L$. We see readily 
that this case can be ignored. Likewise the
case $\delta_1=-1$ with $t\ll K^\varepsilon$ can be dropped from consideration
under Convention 1, because we now have $F/L\gg K^\varepsilon$ and Lemma 6 
works as in the previous case.
\par
On the other hand, if $\delta_1=-1$ and $K^\varepsilon\ll t$, then we may compute
the $v$-integral asymptotically with the saddle point method. The saddle point is
at $v=4t$. We divide the integral into two parts according as $|v-4t|<4t/\rho_1$
and otherwise, with $\rho_1=t^{2/5}$. Then we proceed in a fashion much
similar to the proof of $(3.31)$. Thus, under the current
situation, the integral of $(5.6)$ is seen to be
$$
\eqalignno{
\sim&\sqrt{\pi\over t}\exp\left(-2it\log(4t)+\txt{1\over4}\pi i\right)\cr
&\times\int\xi(f,\ell,w/(4t),4t)\log^2(w/(4t))
\exp\left(iQ_1(4t)\right){dw\over w},&(5.10)
}
$$
where $Q_1(4t)=Q_1\vert_{v=4t}$. 
That is, the estimation of $(5.6)$ has been reduced to that of
$$
{1\over\sqrt{t}}\sum_{\ell\le K^\varepsilon}{1\over\ell}\int
\left|\sum_{f=1}^\infty\phi_2(f){\sigma_1(f)c_\ell(f)\over f^{1-2it}}
\xi(f,\ell,w/(4t),4t)\right|{dw\over w}.\eqno(5.11)
$$
On invoking $(4.42)$, this sum over $f$ is equal to
$$
\eqalignno{
{i\over8\pi^3}&\int_{(0)}\int_{(0)}\!\int_{(2)}
{(4\pi)^{2s_2}(4t)^{s_1-s_2}\over w^{s_1+s_2}
\ell^{2s_2}}\phi^*(s_1)\phi^*(s_2)\phi_2^*(s_3)
C_\ell(s_1-s_2+s_3+1-2it)
\cr
&\times\zeta(s_1-s_2+s_3-2it)\zeta(s_1-s_2+s_3+1-2it)
ds_1ds_2ds_3,&(5.12)
}
$$
where $\phi^*$, $\phi_2^*$ are Mellin transforms of respective functions; and
$$
C_\ell(s)=\sum_{d_1d_2d_3d_4=\ell}{\mu(d_3)\mu(d_4)\over
(d_1d_2d_3^2)^s}d_1d^2_2d_3,\eqno(5.13)
$$
with the M\"obius function $\mu$. In fact, we have, in the region of absolute
convergence,
$$
\eqalignno{
\zeta(s)\zeta(s-1)\sum_{\ell=1}^\infty {C_\ell(s)\over\ell^\lambda}&=
\zeta(s)\zeta(s-1){\zeta(s+\lambda-1)\zeta(s+\lambda-2)
\over\zeta(\lambda)\zeta(2s+\lambda-2)}\cr
&=\sum_{n=1}^\infty{\sigma_1(n)\over
n^s}{\sigma_{1-\lambda}(n)\over\zeta(\lambda)}\cr
&=\sum_{\ell=1}^\infty{1\over\ell^\lambda}
\sum_{n=1}^\infty{\sigma_1(n)c_\ell(n)\over
n^s}, &(5.14)
}
$$
where we have used two well-known formulas of Ramanujan.
Then, after the truncation to $|s_j|\le K^\varepsilon$ for all $j$ in $(5.12)$, we
move the $s_1$-contour to the imaginary axis. Under Convention 1 and $t\gg
K^\varepsilon$, we encounter no singularities. In this way we find that 
$$
S_r\ll K^{\varepsilon}.\eqno(5.15)
$$\medskip
Next, let us consider $S_d$, the contribution of 
the discrete spectrum; note $(2.29)$ and $(2.32)$. We
need first to approximate $\Psi(i\kappa;X)$, where $\kappa\in \B{R}$, and 
$X$ is defined
by
$(4.41)$  with $Y=Q$. To this end, we invoke the identity
$$
\eqalignno{
&{}_2F_1\left(\txt{1\over2}+i\kappa,\txt{1\over2}+i\kappa;
1+2i\kappa;-{1/u}\right)\cr
&=\left(\txt{1\over2}\left(1+\sqrt{1+1/u}\right)\right)^{-1-2i\kappa}
{}_2F_1\left(\txt{1\over2},\txt{1\over2}+i\kappa;1+i\kappa;
\left({1-\sqrt{1+1/u}\over
1+\sqrt{1+1/u}}\right)^2\right),\quad&(5.15)
}
$$
which is an instance of the quadratic transformations of the Gaussian
hypergeometric function (see, e.g., [18, $(9.6.12)$]). This implies that
uniformly in $\kappa$
$$
{}_2F_1 \left(\txt{1\over2}+i\kappa,\txt{1\over2}+i\kappa;1+2i\kappa;-{1/u}\right)
\sim\left(\txt{1\over2}\left(1+\sqrt{1+1/u}\right)\right)^{-1-2i\kappa},
\eqno(5.16)
$$
with $u$ as in $(4.46)$; note that we have currently $u\gg G^2K^{-\varepsilon}$.
Thus, the estimation of $S_d$ is reduced to that of
$$
\eqalignno{
\sum_{f=0}^\infty{\phi_2(f)\over
f^{{1/2}-2it}}&\sum_{\ell=1}^\infty{\phi_1(\ell)c_\ell(f)
\over\ell^{1+2it}}\sum_{j=1}^\infty2^{2\delta_3i\kappa_j}\alpha_j\tau_j(f)H_j^2
\left(\txt{1\over2}\right)\cr
&\times\left(1+{\delta_3i\over\sinh(\pi\kappa_j)}\right)
{\Gamma^2({1\over2}+\delta_3i\kappa_j)\over\Gamma(1+2\delta_3i\kappa_j)}
\Xi(f,\ell,\kappa_j,\delta_1,\delta_3),&(5.17)
}
$$
where 
$$
\eqalignno{
\Xi(f,\ell,\kappa,\delta_1,&\delta_3)=\int\!\int
\xi(f,\ell,u,v)\exp(iQ)\left(\sqrt{u}+\sqrt{u+1}\right)^{-1-2\delta_3i\kappa}
{dudv\over uv^{1+2it}}\qquad &(5.18)
}
$$
with $\delta_3=\pm1$. 
\par
We have
$$
{\partial\over\partial u}\left(Q-2\delta_3\kappa\log
\left(\sqrt{u}+\sqrt{u+1}\right)\right)=2\delta_1{r^2\over u^2v}
-{\delta_3\kappa\over\sqrt{u(u+1)}}.\eqno(5.19)
$$
On the right side, provided $\kappa\gg (K/G)\log K$, 
the second term is dominant, and Lemma 6 becomes relevant with 
$A_0=(M/F)^{-3/2}$, $A_1=M/(F\log^4K)$, $B_1=\kappa F/M$, $\rho\approx M/F$.
In fact, we have in the relevant domain 
$r^2/(u^2v)\ll FK^2L/M^2\ll ((K/G)\log
K)(F/M)$, by $(4.3)$ and $(4.4)$. 
Thus, this case can be ignored. That is, we may
truncate the inner-most sum of
$(5.17)$ to
$\kappa_j\ll (K/G)\log K$ .  We have then
$\left(\sqrt{u}+\sqrt{u+1}\right)^{i\kappa_j}\sim (2\sqrt{u})^{i\kappa_j}$, as
$\kappa_j/u\ll (FK\log K)/(GM)\ll
K^{1+\varepsilon}/G^3\ll K^{-\varepsilon}$ because of $(4.3)$, $(4.40)$,
$(5.3)$, and Convention 1.
\par
Hence the estimation of $S_d$ is further reduced to that of
$$
\sum_{\kappa_j\ll (K/G)\log K}{\alpha_j\over\sqrt{\kappa_j}}
H_j^2\left(\txt{1\over2}\right)\left|
\sum_{f=1}^\infty{\phi_2(f)\tau_j(f)\over f^{{1/2}-2it}}\sum_{\ell=1}^\infty
{\phi_1(\ell)c_\ell(f)\over
\ell^{1+2it}}\Xi_1(f,\ell,\kappa_j,\delta_1,\delta_3)\right|,
\eqno(5.20)
$$
with
$$
\Xi_1(f,\ell,\kappa,\delta_1,\delta_3)=\int\!\int
\xi(f,\ell,u,v)\exp(iQ)
{dudv\over u^{3/2+\delta_3i\kappa}v^{1+2it}}.\eqno(5.21)
$$
As before, we perform the change of variable $u\mapsto w/v$. Then, an application
of the Mellin inversion gives
$$
\eqalignno{
\Xi_1&(f,\ell,\kappa,\delta_1,\delta_3)=-{1\over4\pi^2}
\int_{(0)}\!\int_{(0)}\phi^*(s_1)\phi^*(s_2)
f^{s_2-s_1}(4\pi/\ell)^{2s_2}\cr
&\times\left(\int_0^\infty{\exp(-{1\over2}\delta_1iv)\over
v^{{1/2}-\delta_3i\kappa+2it-s_1+s_2}}dv\right)\left(\int_0^\infty
{\exp(-2\delta_1ir^2/w)\over 
w^{{3/2}+\delta_3i\kappa+s_1+s_2}}dw\right)ds_1ds_2,&(5.22)
}
$$
which can be verified as $(4.25)$. Thus,
$$
\eqalignno{
\qquad\Xi_1(f,\ell,&\kappa,\delta_1,\delta_3)
={2^{-2it}\over8\pi^2 r^{1+2\delta_3i\kappa}}
\int_{(0)}\!\int_{(0)}\phi^*(s_1)\phi^*(s_2)
f^{s_2-s_1}(2\pi/\ell)^{2s_2}\cr
&\times r^{-2s_1-2s_2}\exp\left(\delta_1\pi i\left(
-\delta_3i\kappa+it-s_1\right)\right)\cr
&\times\Gamma\left(\txt{1\over2}+\delta_3i\kappa-2it+s_1-s_2\right)
\Gamma\left(\txt{1\over2}+\delta_3i\kappa +s_1+s_2\right)ds_1ds_2.&(5.23)
}
$$
Inserting this assertion into $(5.20)$, we encounter
$$
\sum_{\ell=1}^\infty{\phi_1(\ell)c_\ell(f)\over\ell^{1+2it+2s_2}}
={1\over2\pi i}\int_{(0)}
\phi_1^*(s_3){\sigma_{-2it-2s_2-s_3}(f)\over\zeta(1+2it+2s_2+s_3)}
ds_3.\eqno(5.24)
$$
\par
Hence 
$$
\eqalignno{
S_d\ll& K^{-1+\varepsilon}
\int_{(0)}\!\int_{(0)}\!\int_{(0)}\left|\phi^*(s_1)\phi^*(s_2)
\phi_1^*(s_3)\right|\cr
\times&\sum_{\kappa_j\ll (K/G)\log K}
{\alpha_j\over\sqrt{\kappa_j}}
H_j^2\left(\txt{1\over2}\right)\left|
\sum_{f=1}^\infty \phi_2(f)\tau_j(f){\sigma_{-2it-2s_2-s_3}(f)\over
f^{{1/2}-2it+s_1-s_2}}\right||ds_1||ds_2||ds_3|\cr
\ll&
K^{-1+\varepsilon}\sup_U\sup_{t_1,t_2}{1\over\sqrt{U}}
\sum_{U\le\kappa_j\le 2U}
\alpha_jH_j^2\left(\txt{1\over2}\right)
\left|\sum_{f=1}^\infty \phi_2(f)\tau_j(f){\sigma_{-it_1}(f)\over
f^{{1/2}-it_2}}\right|,&(5.25)
}
$$
where $U\ll (K/G)\log K$, and $|t_\nu-2t|\ll K^\varepsilon$ $(\nu=1,2)$, after an
obvious truncation of the triple integral. This and Lemma 7 imply that
$$
S_d\ll K^{-1+\varepsilon}(K/G)^{1/2}\left(K/G+F^{1/2}\right)\ll
(K/G^3)^{1/2}K^\varepsilon\ll
K^{\varepsilon},\eqno(5.26)
$$
because $F\ll (K/G)^2 K^\varepsilon$. Here we have used a well-known bound for the
spectral fourth moment of $H_j\left({1\over2}\right)$ that follows from,
e.g., $(1.28)$.
\medskip
The discussion of $S_c$ is analogous to the above, up to $(5.25)$. In fact, the
change is only in that $(5.25)$ is to be replaced by the
expression
$$
K^{-1+\varepsilon}\sup_U\sup_{t_1,t_2}\int_{-U}^{U}
{|\zeta\left({1\over2}+i\kappa\right)|^4\over|\zeta(1+2i\kappa)|^2}
\left|
\sum_{f=1}^\infty \phi_2(f){\sigma_{2i\kappa}(f)\sigma_{-it_1}(f)\over
f^{{1/2}+i\kappa-it_2}}\right|{d\kappa\over
\sqrt{|\kappa|+1}}.\eqno(5.27)
$$
To this we could apply a continuous analogue of Lemma 7, but we
take a different way to motivate a later
purpose. Thus, we note first that the part corresponding to $|\kappa|\ll
K^\varepsilon$ contributes $K^{-1+\varepsilon}F^{1/2}\ll K^\varepsilon/G$. To
treat the part with  $|\kappa|\gg K^\varepsilon$,  we use Mellin inversion. The
last sum is equal to
$$
\eqalignno{
{1\over2\pi i}\int_{(2)}&\phi_2^*(s)\zeta\left(s+\txt{1\over2}+i\kappa-it_2
\right)\zeta\left(s+\txt{1\over2}-i\kappa-it_2\right)
\zeta\left(s+\txt{1\over2}+i\kappa
+i(t_1-t_2)\right)\cr
&\times\zeta\left(s+\txt{1\over2}-i\kappa
+i(t_1-t_2)\right)
\left\{\zeta(s+1+it_1-2it_2)\right\}^{-1}ds,&(5.28)
}
$$
again by a formula of Ramanujan.
After truncating to $|s|\le |\kappa|/2$, we shift the contour to $(0)$.
We see that the integral is
$$
\eqalignno{
\ll &\left\{(1+|\kappa|)^2(1+|\kappa-2t|)
(1+|\kappa+2t|)\right\}^{1/6}K^\varepsilon\cr
+&F^{1/2}\left\{
(1+|\kappa-2t|)^{-1/\varepsilon}+(1+|\kappa+2t|)^{-1/\varepsilon}\right\}
K^{\varepsilon}.&(5.29)
}
$$
The second term comes from the possible simple poles at ${1\over2}\pm i\kappa
+it_2$, since one may assume that $|\kappa|\ge 4|t_1-t_2|$, under
Convention 1.  This implies that the relevant part of the integral in $(5.27)$ is 
$\ll (K/G)^{{7/6}+\varepsilon}$. Hence
$$
S_c\ll K^\varepsilon/G + K^{-1+\varepsilon}(K/G)^{7/6},\eqno(5.30)
$$
which is negligible.
\medskip
It remains to consider $S_h$. From $(2.30)$, it is
$$
\ll\sum_{f=1}^\infty{\phi_2(f)\over f^{1/2}}
\sum_{\ell=1}^\infty{\phi_1(\ell)
\over\ell}|c_\ell(f)|\sum_{\scr{k=6}\atop\scr{2|k}}^\infty
\sum_{j=1}^{\vartheta(k)}\alpha_{j,k}|\tau_{j,k}(f)|H_{j,k}^2
\left(\txt{1\over2}\right)\left|\Psi\left(k-\txt{1\over2};X\right)
\right|.\eqno(5.31)
$$
The $\Psi$-factor is, by $(2.32)$ and $(4.41)$,
$$
\ll{\Gamma(k)^2\over\Gamma(2k)}\int u^{-k-1}
{}_2F_1(k,k;2k;-1/u)du\ll \left({F/M}\right)^{k-\varepsilon}.\eqno(5.32)
$$
Here the range of integration is given by
$(4.46)$, and the bound is uniform in $k$; the latter can be seen by using
Gauss' integral representation of ${}_2F_1$. Also, we
invoke Deligne's bound
$$
|\tau_{j,k}(f)|\le d(f),\eqno (5.33)
$$
and its elementary consequence
$$
H_{j,k}\left(\txt{1\over2}\right)\ll 
k^{{1/2}+\varepsilon}.\eqno(5.34)
$$
Further,
$$
\sum_{j=1}^{\vartheta(k)}\alpha_{j,k}\ll k\eqno(5.35)
$$
(see [23, Lemma 3.3]).
One could replace $(5.33)$ by the bound given in [23, $(3.1.22)$]
for instance, and $(5.34)$ by an easier convexity bound. At any event, the
above combination gives that
$$
S_h\ll F^{1/2}L(F/M)^6K^\varepsilon,\eqno(5.36)
$$
which is negligible.
\medskip
Collecting $(4.45)$, $(5.1)$, $(5.15)$, $(5.26)$, $(5.30)$, and $(5.36)$, we
end the proof of $(1.27)$ on the condition $(5.3)$. In particular, we have
proved $(1.23)$ and consequentially Ivi\'c's bound for
$H_j\left({1\over2}\right)$ as well, in a wider context. 
\vskip 1cm
\noindent
\centerline{\bf 6. Intermediate range}
\bigskip
\noindent
Now, we enter into the intermediate range; or more precisely we shall consider
${\cal S}_4$, with $(5.1)$, under the conditions
$$
\hbox{$K^{1+\varepsilon}/G\ll t\ll K^{{3/2}-\varepsilon}\,$ and
$\,|K-t|\gg G$}.\eqno (6.1)
$$
Here the quantity $K^{1+\varepsilon}/G$ should be equal to the
same in $(5.3)$, because of an obvious reason. The second condition is by no means
a restriction, because $(3.50)$ already gives what we desire. 
The conditions $(1.26)$, $(4.3)$, $(4.4)$ with $(3.38)$, and $(4.5)$ are of course
retained, but $(4.40)$ now becomes 
$$
F\ll tLK^\varepsilon.\eqno(6.2)
$$
We shall be brief occasionally, since the reasoning is 
analogous, though not quite, to that developed in the preceding section.
\medskip
We begin with $S_r$. We consider, instead of $(5.6)$, the expression
$$
\sum_{f=1}^\infty{\phi_2(f)\sigma_1(f)\over f^{1-2it}}
\sum_{\ell=1}^\infty{\phi_1(\ell)c_\ell(f)\over \ell^{1+2it}}
\int\!\int\xi(f,\ell,w/v,v){(\log w/v)^2\over 
v^{1+2it}w}\exp\left(Y_1i\right)
dv\,dw,\eqno(6.3)
$$
where the integration range is $(5.7)$ and
$$
Y_1=-{tv\over2w}\left(1-{v\over4 w}\right)-{1\over2}\delta_1v
\left(1+{v^2\over32w^2}\right)-2\delta_1{r^2\over w}
\left(1-{v\over4w}\right).\eqno(6.4)
$$
We integrate first with respect to $v$. We have 
$$
{\partial\over\partial v}(Y_1-2t\log v)=-2{t\over v}\left(1
+{v\over4w}-{v^2\over 8w^2}\right)-{1\over2}\delta_1
\left(1+{3\over32}{v^2\over w^2}-{r^2\over w^2}\right)
\eqno(6.5)
$$
which is close to $-2t/v-{1\over2}\delta_1$;
in fact, $v/w\ll F/M\ll tK^\varepsilon/(GK)\ll t^{-1/3}$, and $r^2/w^2\ll
(G/\log K)^{-2}$. Hence, if $\delta_1=1$, then $(6.3)$ is negligible, as
can be confirmed with the procedure following $(5.9)$.
Thus, let us assume that $\delta_1=-1$.
Then, we apply the saddle point method. 
The saddle point is at $v_0\sim 4t$ or more precisely it satisfies the recursive
equation
$$
v_0=4t\left(1
+{v_0\over4w}-{v_0^2\over 8w^2}\right)
\left(1+{3\over32}{v_0^2\over w^2}-{r^2\over w^2}\right)^{-1}.\eqno(6.6)
$$
When $|v-v_0|<v_0/\rho_1$, with $\rho_1=t^{{2/5}}$, we have
$$
\eqalignno{
Y_1-2t\log v&=Y_1(v_0)-2t\log v_0
+2t\sum_{j=2}^\infty{1\over
j}\left({v_0-v\over v_0}\right)^j\cr
&+{1\over2}Y_1^{(2)}(v_0)(v-v_0)^2+{1\over6}Y_1^{(3)}(v_0)(v-v_0)^3,&(6.7)
}
$$
where $Y_1^{(\nu)}(v_0)=((\partial/\partial v)^\nu)_{v=v_0} Y_1$. We note that
$Y_1^{(2)}(v_0)\ll t/w^2\ll tK^\varepsilon/(GK)^2$, and 
$Y_1^{(3)}(v_0)\ll K^\varepsilon/(GK)^2$.
Thus, 
$$
Y_1-2t\log v\sim Y_1(v_0)-2t\log v_0
+t\left({v_0-v\over v_0}\right)^2,\eqno(6.8)
$$
and the integral of $(6.3)$ is 
$$
\sim e^{{1\over4}\pi
i}\sqrt{\pi\over t}\int\xi(f,\ell,w/v_0,v_0)\log^2(w/v_0)
\exp\left(iY_1(v_0)
-2it\log v_0\right){dw\over w}.\eqno(6.9)
$$ 
This corresponds to $(5.10)$.
\par
We shall show that $(6.9)$ is negligibly small, if 
$\ell\gg K^\varepsilon$. To this end, we note that
$$
{d\over dw}\left(Y_1(v_0)-2t\log v_0\right)=
{tv_0\over2w^2}\left(1-{v_0\over2 w}\right)-
{v_0^3\over32w^3}-2{r^2\over w^2}
\left(1-{v_0\over 2w}\right),\eqno(6.10)
$$
since the left side is equal to $\left(Y_1\right)_w(v_0)$ by the definition of
$v_0$. Inserting the approximation $v_0=4t\left(1+{t/ w}\right)
+O\left({t(t+K)^2/ w^2}\right)$, which follows readily from $(6.6)$, we get
$$
{d\over dw}\left(Y_1(v_0)-2t\log v_0\right)=2{t^2-r^2\over w^2}\left(1-2{t\over
w}\right)+O\left({(t(t+K))^2\over w^4}\right).\eqno(6.11)
$$
Thus,
$$
{d\over dw}\left(Y_1(v_0)-2t\log v_0\right)
\approx {T\over w^2}.\eqno(6.12)
$$
In fact, $t^2-r^2\approx T$, and 
$$
\eqalignno{
{(t(t+K))^2\over w^2T}&\ll {t^2(t+K)\over(GK)^2|K-t|}\cr
&\ll{t^2\over(t+K)^{5/3}|K-t|}K^{-\varepsilon}
\ll {t^{1/3}\over|K-t|}K^{-\varepsilon},&(6.13)
}
$$
which is negligibly small; here we have used $(1.26)$ and the second condition
in $(6.1)$. Then, Lemma 6 is applied to
$(6.9)$, with $A_0=(GK)^{-1}(\log K)^2$, $A_1=M/(L(\log K)^4)$, 
$B_1=T(L/M)^2$, $\rho\approx M/L$. The integral is $\ll (LT/M)^{-P}(\log K)^2$.
Hence, $(6.9)$ is indeed negligibly small if $L\gg K^\varepsilon$, because we have
$(4.3)$.
\par
That is, in dealing with $(6.3)$, we may assume that we have $L\ll K^\varepsilon$
together with $(6.9)$. We are, thus, in a situation much analogous to that with
$S_r$ in the lower range. In this way we are led again to
$$
S_r\ll K^\varepsilon.\eqno(6.14)
$$
\medskip
We turn to $S_d$. The reduction to $(5.17)$--$(5.18)$ does not need to be altered,
except for the replacement of $Q$ by $Y$. With this, we shall consider $\Xi$.
We have
$$
\eqalignno{
&{\partial\over\partial u}\left(Y-2\delta_3\kappa\log(\sqrt{u}
+\sqrt{u+1})\right)\cr
=&{t\over 2u^2}\left(1-{1\over 2u}\right)+{1\over32}\delta_1
{v\over u^3}+2\delta_1{r^2\over u^2v}\left(1-{1\over 2u}\right)
-{\delta_3\kappa\over\sqrt{u(u+1)}}.&(6.15)
}
$$
This shows in particular that the part with $\kappa_j\gg K^{1/\varepsilon}$ of
$S_d$ can be discarded, as Lemma 6 implies. Thus we have
the initial truncation $\kappa_j\ll K^{1/\varepsilon}$. 
\par
We then integrate with respect to $v$. We have
$$
{\partial\over\partial v}\left(Y-2t\log v\right)=-{1\over2}\delta_1
\left(1+{1\over32u^2}\right)+2\delta_1{r^2\over uv^2}\left(1-{1\over4u}
\right)-{2t\over v}.\eqno(6.16)
$$
Note that $r^2/(uv)\ll K^2L/M\ll
(K/G)\log K\ll tK^{-\varepsilon}$ because of $(4.3)$ and the first condition in
$(6.1)$. Hence, we may adopt the argument following
$(5.9)$, and see that the case $\delta_1=1$ can be discarded. 
Hereafter we shall assume that $\delta_1=-1$. The $v$-integral has a
saddle point at $v_1$, which satisfies the recursive equation
$$
v_1=4t\left(1+{r^2\over uv_1t}\left(1-{1\over4u}\right)\right)
\left(1+{1\over
32u^2}\right)^{-1};\eqno(6.17)
$$
and $v_1=4t\left(1+O\left(K^{-\varepsilon}\right)\right)$; in particular, $(5.7)$
gives
$$
F\approx tL,\eqno(6.18)
$$
which replaces $(6.2)$. The saddle point 
method yields that $\Xi$ defined by $(5.18)$ with $Q=Y$ is
$$
\sim e^{{1\over4}\pi
i}\sqrt{\pi\over t}\int\xi(f,\ell,u,v_1)
\exp\left(iY(v_1)
-2it\log v_1\right)\left(\sqrt{u}+\sqrt{u+1}
\right)^{-1-2\delta_3i\kappa}{du\over
u},\eqno(6.19)
$$
where $Y(v_1)=Y|_{v=v_1}$.
In this, we have, by the definition of $v_1$,
$$
{\partial\over\partial u}\left(Y(v_1)-2t\log v_1-2\delta_3\kappa
\log\left(\sqrt{u}+\sqrt{u+1}\right)\right)=\left(Y\right)_u(v_1)
-{2\delta_3\kappa\over\sqrt{u(u+1)}}.\eqno(6.20)
$$
We shall show that
$$
\left(Y\right)_u(v_1)={t\over 2u^2}\left(1-\left({r\over t}\right)^2\right)
\left(1+O(K^{-\varepsilon})\right)\approx {T\over u^2t}.\eqno(6.21)
$$
It suffices to prove the asymptotics; and to this end
we may assume that $t\ge K^{1-\varepsilon}$, since otherwise the assertion follows
immediately from the first line of $(6.22)$ below. Then, we note that $(6.17)$
gives $v_1/(16ut)=1/(4u)+O(K^\varepsilon/u^2)$ and
$4r^2/(v_1t)=(r/t)^2(1-(r/t)^2/(4u))+O(K^\varepsilon/u^2)$.
Thus,
$$
\eqalignno{
\left(Y\right)_u(v_1)=&{t\over 2u^2}\left(1-{1\over2u}-{v_1\over16ut}
-4{r^2\over v_1t}\left(1-{1\over2u}\right)\right)\cr
=&{t\over 2u^2}\left(1-\left({r\over
t}\right)^2\right)\left(1-{3\over4u}-{1\over4u}
\left({r\over t}\right)^2\right)+O\left(tK^\varepsilon/u^4\right).&(6.22)
}
$$
Also, by $(1.26)$ and Convention 1,
$$
{1\over u^2}\left(1-\left({r\over t}\right)^2\right)^{-1}\ll {t^{1/3}\over
|K-t|}K^{-\varepsilon},\eqno(6.23)
$$
which proves $(6.21)$, because of the second condition in $(6.1)$. In passing, we
stress that both assumptions in $(6.1)$ are indeed required in the above.
\par
With $(6.20)$--$(6.21)$, Lemma 6 allows us to impose the truncation
$$
\kappa_j\ll \kappa_0=TK^\varepsilon/(GK).\eqno(6.24)
$$
In fact, we may set $\rho\approx M/F$, and in the relevant domain of $u$ we have, 
$T/(ut)\ll TK^\varepsilon/(GK)$. Thus, provided $\kappa\gg\kappa_0$, a
specification is given by $A_0=(M/F)^{-3/2}$, $A_1=M/(F\log^4K)$, 
$B_1=\kappa F/M$. Note that $\kappa_0\gg K^\varepsilon$, under Convention 1.
\par
Hence, the estimation of $S_d$ has been reduced to that of
$$
{1\over\sqrt{tU}}\sum_{U\le\kappa_j\le2U}\alpha_j
H_j^2\left(\txt{1\over2}\right)\left|
\sum_{f=1}^\infty{\phi_2(f)\tau_j(f)\over f^{{1/2}-2it}}\sum_{\ell=1}^\infty
{\phi_1(\ell)c_\ell(f)\over
\ell^{1+2it}}\Xi_2(f,\ell,\kappa_j,\delta_3)\right|,
\eqno(6.25)
$$
where $U\ll \kappa_0$, and
$$
\eqalignno{
\Xi_2(&f,\ell,\kappa,\delta_3)=\int_0^\infty \varkappa(u)
\xi(f,\ell,u,v_1)\cr
&\times\exp(iY(v_1)-2it\log v_1)
\left(\sqrt{u}+\sqrt{u+1}\right)^{-1-2i\delta_3\kappa}{du\over u}.&(6.26)
}
$$
Here $\varkappa$ is a smooth weight whose r\^ole is analogous to that of $\theta$
in $(3.26)$.
\par
Appealing to Mellin inversion, we find that $(6.25)$ is
$$
\eqalignno{
&\ll{K^\varepsilon\over\sqrt{tU}}\int_{(0)}\!\int_{(0)}\!\int_{(0)} 
\left|\phi^*(s_1)\phi^*(s_2)\phi_1^*(s_3)\right|\sum_{U\le\kappa_j\le 2U}
\alpha_jH_j^2\left(\txt{1\over2}\right)
\cr
&\times
\left|\Lambda(\kappa_j,\delta_3;s_1,s_2)\right|\left|\sum_{f=1}^\infty
\phi_2(f)\tau_j(f){\sigma_{-2it-2s_2-s_3}(f)\over
f^{{1/2}-2it+s_1-s_2}}\right||ds_1||ds_2||ds_3|,&(6.27) }
$$
where
$$
\eqalignno{
\Lambda(\kappa,\delta_3;s_1,s_2)=\int_0^\infty&\varkappa(u)
\exp(iY(v_1)-2it\log v_1)\cr
&\times\left(\sqrt{u}+\sqrt{u+1}\right)^{-1-2i\delta_3\kappa}{du\over
u^{1+s_1+s_2}}. &(6.28)
}
$$
This should be compared with $(5.25)$.
\medskip
We shall bound $\Lambda$. Naturally we could truncate $(6.27)$ to
$|s_1|,\,|s_2|,\,|s_3|\ll K^\varepsilon$, as we assume now. We have, by
$(6.20)$--$(6.21)$,
$$
\eqalignno{
{\partial\over\partial u}&\left\{Y(v_1)-2t\log v_1+i(s_1+s_2)\log u-2\delta_3\kappa
\log(\sqrt{u}+\sqrt{u+1})\right\}\cr
&=2\pi^2\sgn(t-r){T\over tu^2}\left(1+O(K^{-\varepsilon})
\right)+i{s_1+s_2\over u}-{2\delta_3\kappa
\over\sqrt{u(u+1)}}.&(6.29)
}
$$
When $\delta_3=\sgn(r-t)$ and $L+\kappa\gg K^\varepsilon$, we may appeal to
Lemma 6 with $A_0=(M/F)^{-3/2}$, $A_1=M/F$, $B_1=(TL/M+\kappa)F/M$, $\rho\approx
M/F$. In fact, this assertion on $B_1$ follows from the Taylor expansion of the
left side of $(6.29)$ around any real point in the relevant domain, coupled with
the fact that $ut\approx uv_1\approx M/L$ by $(4.46)$, and $T\gg M$ as well as
$|s_1+s_2|\ll K^\varepsilon$, thus under Convention 1. We find that $\Lambda$ is
negligibly small with the present supposition. That is,  provided
$\delta_3=\sgn(r-t)$, we may impose the truncation $L+\kappa_j\ll K^\varepsilon$
in $(6.27)$; in particular,
$U\ll K^\varepsilon$, and $F\ll tK^\varepsilon$ by
$(6.18)$. Then, applying Lemma 7 to the sum over $\kappa_j$ of $(6.27)$ we
immediately find that the case $\delta_3=\sgn(r-t)$ can be dropped.
\par
Hence we assume now that $\delta_3=\sgn(t-r)$. With this, let $u_1$ be the saddle
point of the integral in $(6.28)$. By the second line of $(6.29)$ we have $1/u_1\ll
(t/T) (\kappa+K^\varepsilon)$. This implies that if
$\kappa_j\ll K^\varepsilon$, then $L\ll K^\varepsilon$, because $(6.18)$
gives $tL/M\approx F/M\approx 1/u_1$. That is, we can ignore this situation as
well; and hence we may assume that $\delta_3=\sgn(t-r)$ and $\kappa\gg
K^\varepsilon$. Then, the saddle point method yields, in a fashion similar
to the argument leading up to $(3.32)$, that
$$
\Lambda(\kappa,\sgn(t-r);s_1,s_2)\ll{1\over\sqrt{u_1^3|\Lambda_0|}},\eqno(6.30)
$$
where
$$
\Lambda_0=\left({\partial\over\partial u}\right)^2_{u=u_1}
\left\{Y(v_1)-2t\log v_1+i(s_1+s_2)\log u+2\kappa
\log(\sqrt{u}+\sqrt{u+1})\right\},\eqno(6.31)
$$
with $\delta_1=\sgn(t-r)$ in the definition of $Y$. We have
$$
\Lambda_0\approx {T\over tu_1^3}.\eqno(6.32)
$$
Inserting this into $(6.27)$ via $(6.30)$ together with the
aforementioned truncation of the triple integral, we find that
$$
S_d\ll K^\varepsilon\sup_U\sup_{t_1,t_2}{1\over\sqrt{TU}}
\sum_{U\le\kappa_j\le 2U}
\alpha_jH_j^2\left(\txt{1\over2}\right)
\left|\sum_{f=1}^\infty \phi_2(f)\tau_j(f){\sigma_{-it_1}(f)\over
f^{{1/2}-it_2}}\right|,\eqno(6.33)
$$
with 
$$
U\ll  TK^\varepsilon/(GK),\quad|t_\nu-2t|\ll K^\varepsilon\;(\nu=1,2),\quad
F\ll Ut.\eqno(6.34)
$$
The bound for $F$ follows from the observation that $M/F\approx u_1\approx
T/(t\kappa)$. The assertion $(6.33)$ is obviously an extension of $(5.25)$,
since $T\approx K^2$ in the lower range.
\medskip
The discussion of $D_c$ and $D_h$ is analogous to that in the previous section; and
it can readily be seen that their contributions are again negligible.
\medskip
Collecting $(6.14)$, $(6.33)$, and the last assertion, we conclude that our problem
has been reduced to the estimation of the spectral sum in $(6.33)$. However,
unlike the case of the lower range, the sole application of Lemma 7 to $(6.33)$
does not settle our problem. In fact, we end up with
$$
{\cal S}(G,K)\ll \left(GK+\sqrt{Tt}\right)^{1+\varepsilon},\eqno(6.35)
$$
which yields the inferior exponent $3\over8$ in place of $1\over3$ in $(1.1)$.
\medskip
To resolve this difficulty, we
have to devise yet another spectral mean value result, whose discussion
is to be developed in the next section. To make our next aim
clearer, we perform a transformation analogous to $(5.28)$ to the sum over $f$ in
$(6.33)$. It is expressed as
$$
{1\over 2\pi i}\int_{(0)}\phi_2^*(s)
{H_j\left(s+{1\over2}-it_2\right)
H_j\left(s+{1\over2}+i(t_1-t_2)\right)
\over\zeta(2s+1+it_1-2it_2)}
ds,\eqno(6.36)
$$
which can be truncated to $|s|\ll K^\varepsilon$. 
Hence, in the intermediate range we have 
$$
{\cal S}(G,K)\ll GK^{1+\varepsilon}
\left(1+\sup_U\sup_{t_3}\sqrt{{\cal T}(U,t_3)(U/T)}\right),\eqno(6.37)
$$
where $t_3\approx t$, and $U\ll TK^\varepsilon/(GK)$, and ${\cal T}$ is
defined by $(1.29)$. Appealing to Theorem 2, we would be able to end
the proof of $(1. 27)$ immediately. 
\vskip 1cm
\noindent
\centerline{\bf 7. Hybrid moment}
\bigskip
\noindent
Now, we begin our discussion of the mean value ${\cal T}(K,t)$. Our aim
is to prove $(1.30)$. In the course of discussion we
shall encounter two instances of applications of Lemma 4, as indicated in
Introduction. Accordingly, the present section is divided into two parts,
with the second starting at $(7.27)$. It should be understood 
that the basic parameters are independent of 
those utilised in the above. On the other hand,
smooth weights attached to sums over integers are as before, and bounds for their
derivatives will be applied without mention. Also, we shall not give details
about applications of Lemma 6 and the saddle point method, since they are 
much similar to those we have encountered in the above.
\medskip
First of all, we observe that we may restrict ourselves to the situation
$$
K^{1+\varepsilon}\ll t\ll K^{2-\varepsilon},\eqno(7.1)
$$
with a sufficiently large $K$. In fact, the case 
$0\le t\ll K^{1+\varepsilon}$ is contained in
$(1.28)$. On the other hand, the case $K^{2-\varepsilon}\ll t$ is readily 
settled with a combination of Ivi\'c's bound for $H_j\left({1\over2}\right)$ and
$(1.19)$. Thus $(7.1)$ is assumed throughout the present section.
\medskip
Then, we put
$$
\eqalignno{
h(r)&=K^{-2P_0}\prod_{p=0}^{P_0-1} 
\left(r^2+\left(p+\txt{1\over2}\right)^2\right)\cr
&\times
\left\{\exp\left(-\left(({r-K)/G}\right)^2\right)+
\exp\left(-\left({(r+K)/ G}\right)^2\right)\right\},&(7.2)
}
$$
with  
$$
G= K^{1-\varepsilon},\eqno(7.3)
$$ 
and consider 
$$
\sum_{j=1}^\infty\alpha_jH_j^2\left(\txt{1\over2}\right)
\left|H_j\left(\txt{1\over2}+it\right)\right|^2h(\kappa_j).\eqno(7.4)
$$
Obviously it suffices to prove that this is bounded by the right side of
$(1.30)$. The integer $P_0\ge1$ is to be fixed later, at $(7.30)$, where the effect
of the polynomial factor  of high order will become apparent.
\par
In much the same way as $(3.40)$ one may show that 
$$
H_j\left(\txt{1\over2}+it\right)
\ll (\log K)^2\sum_{M\ll t}
\int_{\gamma^{-1}-i\gamma^2}^{\gamma^{-1}+i\gamma^2}
\left|\sum_{m=1}^\infty
\phi(m;M)\tau_j(m)m^{-{1/2}-it-\xi}\right|
|d\xi|,\eqno(7.5)
$$
with $\gamma=\log^2K$ and dyadic numbers $M$, uniformly in $\psi_j$ under
consideration. Thus, in place of the factor
$\left|H_j\left({1\over2}+it\right)\right|^2$ of $(7.4)$, we may put
$$
\sum_{m=1}^\infty\sum_{n=1}^\infty{\phi_0(m)
\overline{\phi_0(n)}\over
\sqrt{mn}}\left({m/n}\right)^{it}\tau_j(m)\tau_j(n),\eqno(7.6)
$$ 
where $\phi_0$ is as in $(4.1)$ with $q=1$. We apply Mellin inversion to $\phi_0$
and subsequent truncation of the integration range, and also invoke
$(1.6)$. In this way, our problem $(7.1)$ is reduced to the estimation of 
$$
\sum_{n=1}^\infty\phi(n){\sigma_{2it_1}(n)\over n^{{1/2}+it_2}}
\sum_{j=1}^\infty\alpha_j\tau_j(n)H_j^2\left(\txt{1\over2}\right)
h(\kappa_j),\eqno(7.7)
$$
where $|t_\nu-t|\ll K^\varepsilon$ $(\nu=1,2)$, and $\phi$ is as in $(4.2)$ but
with $M\ll t^2$.
\medskip
We apply Lemma 3 or $(2.9)$ to the last inner sum. We have
$$
\eqalignno{
&\e{H}_1(n;h)\ll d(n)n^{-{1/2}}GK\log^2\!K,\quad \e{H}_3(n;h)
\ll\exp(-\log^2\!K),\cr
&\e{H}_5(n;h)\ll d(n)n^{-{1/2}}\exp(-K),\quad
\e{H}_6(n;h)\ll\sigma_{-1}(n)n^{1/2}\exp\left(-(K/G)^2\right),&(7.8)
}
$$
which can be verified following the discussion on [23, p.\ 128]. As to the
contribution of $\e{H}_7$, it is expressed as
$$
\eqalignno{
&{i\over2\pi^2}\int_{-\infty}^\infty h(r){|\zeta(\txt{1\over2}+ir)|^4\over
|\zeta(1+2ir)|^2}\int_{(1)}\phi^*(s)\zeta\left(s+\txt{1\over2}-ir+it_2\right)
\cr
&\times\zeta\left(s+\txt{1\over2}+ir+it_2\right)
\zeta\left(s+\txt{1\over2}-ir+it_2-2it_1\right)\cr
&\times\zeta\left(s+\txt{1\over2}+ir+it_2-2it_1\right)
\left\{\zeta\left(2s+1+2it_2-2it_1\right)\right\}^{-1}ds\,dr,&(7.9)
}
$$
with the Mellin transform $\phi^*$ of $\phi$. We truncate the inner integral
to $s\ll K^\varepsilon$, and shift the contour to the imaginary axis. Because
of the lower bound for $t$ in $(7.1)$, we do not encounter any sigularity, and
$(7.9)$ is seen to be $\ll Kt^{2/3}$.
Hence, the contribution to $(7.7)$ of the terms in $(2.9)$ except 
for $\e{H}_2$ and $\e{H}_4$ is  $\ll GK\log^5\!K+Kt^{2/3}$.
\medskip
The discussion of the contribution of $\e{H}_2$ and $\e{H}_4$ is
remaining. We shall treat $\e{H}_2$ first. 
The argument is an adaptation 
of  [12, Chapters 2--3] and [23, Section 3.4], except for the treatment 
of a contribution of holomorphic cusp
forms, i.e., $(7.22)$ below. 
\medskip
Thus, we need to estimate the sum
$$
\sum_{f=1}^\infty {d(f)\over f^{1+it_2}}\sum_{n=1}^\infty
\sigma_{2it_1}(n)d(n+f)W^+_f(n/f),\eqno(7.10)
$$
where 
$$
W^+_f(x)=\phi(fx)x^{-{1/2}-it_2}\Psi^+(1/x;h),\eqno(7.11)
$$
with $\Psi^+$ defined by $(2.7)$. 
We may truncate the outer sum to $f\ll K^{1/\varepsilon}$, as can be
confirmed by moving the $r$-contour of 
$(2.7)$ to $\Im(r)=-{3\over4}$. We have, with $x=n/f$, 
$$
\eqalignno{
\Psi^+(1/x;h)\sim &4\pi^{3/2}
GK\,\Re\int_0^1\left(y(1-y)(1+xy)\right)^{-{1/2}}\cr
&\times\left({y(1-y)\over 1/x+y}\right)^{iK}
\exp\left(-\left({G\over2}\log{y(1-y)\over
1/x+y}\right)^2\right)dy.&(7.12)
}
$$
The maximum value of 
$y(1-y)/(1/x+y)$ is less than $\left(1+x^{-1/2}\right)^{-2}$.
Thus, $\Psi^+(1/x;h)$ is negligibly small when $x\ll G^2K^{-\varepsilon}$.
Therefore we may further impose the truncation $f\ll K^\varepsilon M/G^2$.
With this, we compute the last integral asymptotically by the
saddle point method. We get, after some simplification,
$$
\Psi^+(1/x;h)\sim 
2^{3/2}\pi^2GK^{1/2}x^{-{1/4}}\Re\exp(2iK/\sqrt{x}-G^2/x).\eqno(7.13)
$$
Then, via $(7.11)$, we may consider, instead of $(7.10)$, the estimation of
$$
GK^{1/2}\sum_{f=1}^\infty {\phi_2(f)d(f)\over f^{1+it_2}}\sum_{n=1}^\infty
\sigma_{2it_1}(n)d(n+f)V^+_f(n/f,\delta),\eqno(7.14)
$$
where 
$$
V^+_f(x,\delta)=\phi(fx)x^{-{3/4}-it_2}\exp(2i\delta K/\sqrt{x}
-G^2/x),\quad \delta=\pm1,\eqno(7.15)
$$
and $\phi_2$ is as in $(4.39)$ but with 
$$
F\ll K^\varepsilon M/G^2.\eqno(7.16)
$$
\medskip
To the inner sum of $(7.14)$ we apply the spectral decomposition $(2.17)$, with
$\alpha=2it_1$ and $\beta=0$. The leading term $D_r$ has the factor
$Y_f(x;2it_1,0)$ in the integrand, and it should be understood as a limit of
the right side of $(2.22)$. The definition $(7.15)$ implies that Lemma 6
applies here and that the contribution of $D_r$ is negligible.
\medskip
We turn to the contribution of the part $D_d$ defined by $(2.19)$. 
Thus, let us study the $\Psi_\pm$-factors. By $(2.24)$--$(2.25)$, 
$$
\eqalignno{
&\Psi_+\left(i\kappa;2it_1,0;V^+_f(\cdot;\delta)\right)={1\over4\pi
i}\cosh(\pi t_1)\int_{(\varepsilon)}
\cos(\pi s)\Gamma(s+i\kappa)\Gamma(s-i\kappa)\cr
&\times\Gamma\left(\txt{1\over2}-it_1-s\right)
\Gamma\left(\txt{1\over2}+it_1-s\right)
\left\{\int_0^\infty x^{s+it_1-{1/2}}V^+_f(\cdot;\delta)dx\right\}ds,
&(7.17)
 }
$$
and
$$
\eqalignno{
&\Psi_-\left(i\kappa;2it_1,0;V^+_f(\cdot;\delta)\right)={1\over4\pi
i}\cosh(\pi\kappa)\int_{(\varepsilon)}
\sin(\pi s)\Gamma(s+i\kappa)\Gamma(s-i\kappa)\cr
&\times\Gamma\left(\txt{1\over2}-it_1-s\right)
\Gamma\left(\txt{1\over2}+it_1-s\right)
\left\{\int_0^\infty x^{s+it_1-{1/2}}V^+_f(\cdot;\delta)dx\right\}ds,
&(7.18)
 }
$$
with $\kappa$ being real. We shift the $s$-contour to
$(-P_1+{1\over4})$ with $P_1$ a non-negative integer. On the new contour the
integrands are 
$$
\eqalignno{\ll
&\exp\Big(\pi|\lambda|-\pi\max(|\lambda|,|\kappa|)-\pi\max(|\lambda|,t_1)\Big)
\cr
\times&\left(1+|\lambda|K^{-\varepsilon}\right)^{-P_2}\left({F\over
M}\cdot{(1+|\lambda+t_1|)(1+|\lambda-t_1|)\over
(1+|\lambda+\kappa|)(1+|\lambda-\kappa|)}\right)^{P_1},&(7.19)
}
$$
with $\Im(s)=\lambda$ and any integer $P_2\ge0$, 
where the implied constant depends only on 
$\varepsilon$, $P_1$, and $P_2$; the residues are bounded analogously, with
$|\lambda|=|\kappa|$ and $j< P_1$ in place of the exponent $P_1$. In
fact, the Gamma-factors are easy  to bound, and the $x$-integral is
$\ll \left(1+|\lambda|K^{-\varepsilon}\right)^{-P_2}(F/M)^{P_1}$, 
which can be confirmed by Lemma 6, while noting
$(7.3)$, $|t_1-t_2|\ll K^\varepsilon$, and $x\approx M/F$ with $(7.16)$. Thus,
taking $P_1$ sufficiently large, we see that the truncation $\kappa_j\ll
K^{1/\varepsilon}$ can be introduced. With this, we set $P_1=0$, getting the
truncation $\lambda\ll K^\varepsilon$  and subsequently $\kappa_j\ll
K^\varepsilon$, under Convention 1.  Then we shift the truncated $s$-contour to
$\left(P_3\right)$ with an integer $P_3\ge0$.  No singularities are encountered.
On the new contour, the
$s$-integrands are bounded by $(7.19)$, but with the exponent $P_1$ being replaced
by $-P_3$. We get
$$
\Psi_\pm(i\kappa;2it_1,0;V^+_f)\ll \left({M\over F}
\cdot{K^\varepsilon\over t^2}\right)^{P_3},\quad 
|\kappa|\ll K^\varepsilon.\eqno(7.20)
$$
Taking $P_3$ sufficiently large, we see that we may impose the
truncation $F\ll K^\varepsilon$ as well.
Then, we insert $(7.15)$ into $(7.17)$--$(7.18)$, and move the $f$-sum
innermost.  The contribution of $D_d$ to $(7.13)$ is
now seen to be
$$
\eqalignno{
\ll GK^{1/2+\varepsilon}&\sum_{\kappa_j\ll
K^\varepsilon}\alpha_j\left|H_j\left(\txt{1\over2}
+it_1\right)\right|^2\cr
&\times\int_{t^2/K^{\varepsilon}}^\infty
\left|\sum_{f\ll K^\varepsilon}\phi(fx)\phi_2(f){d(f)\tau_j(f)\over
f^{{1/2}+i(t_2-t_1)}}\right|{dx\over x^{{5/4}}},
&(7.21)}
$$
where the lower bound for $x$ is due to
$(7.20)$.  Then, we appeal to either Meurman's bound or $(1.18)$ for Hecke
$L$-functions as well as to a uniform bound for $\tau_j(f)$ (see [23,
$(3.1.18)$]), and  find that $(7.21)$ is $\ll  K^{3/2}t^{1/6}$, which implies that
the $D_d$ part of
$(7.14)$ is negligible.
\medskip
The discussion of the contribution of $D_c$ defined by $(2.21)$ is analogous to the
above, and can be skipped. It should be noted that the assertion $(7.20)$, i.e.,
$F\ll K^\varepsilon$, is
actually not necessary when we deal with the part involving $D_d$. The truncation
$\kappa_j\ll K^\varepsilon$ suffices. However, the part involving $D_c$ requires
$(7.20)$, for the function $f^{-i\kappa}\sigma_{2i\kappa }(f)$ with a small
$\kappa$ is non-oscillating, unlike $\tau_j(f)$. 
\medskip
As to the $D_h$ part of $(7.14)$, we need to bound
$$
\eqalignno{
GK^{1/2}\sum_{f=1}^\infty{\phi_2(f)d(f)\over f^{{1/2}
+i(t_2-t_1)}}
&\sum_{k=6}^\infty\sum_{j=1}^{\vartheta(k)}(-1)^k\alpha_{j,k}\tau_{j,k}(f)
\left|H_{j,k}\left(\txt{1\over2}+it_1\right)\right|^2\cr
&\times\Psi_+\left(k-\txt{1\over2};2it_1,
0;V_f^+(\cdot;\delta)\right).&(7.22)
}
$$ 
Apart from a constant multiplier, the factor $\Psi_+$ is equal to 
$$
\eqalignno{
\cosh&(\pi t_1)\int_{(0)}{\Gamma(k+s-{1\over2})\over\Gamma(k+{1\over2}-s)}
\Gamma\left(\txt{1\over2}+it_1-s\right)
\Gamma\left(\txt{1\over2}-it_1-s\right)\cr
&\times\int_0^\infty \phi(fx)
x^{s-{5/4}+i(t_1-t_2)}
\exp(2i\delta K/\sqrt{x}-G^2/x)dxds.&(7.23)
}
$$
We suppose first that $k\gg K^{1/\varepsilon}$, and shift
the $s$-contour far to the left, without passing over any pole.  We see readily
that the integral is negligibly small. Thus, we get the truncation of
$(7.22)$ to $k\ll K^{1/\varepsilon}$. This implies that we may introduce the
truncation $|s|\ll K^\varepsilon$ as before. Then, shifting the 
truncated $s$-contour far to the right, we see that
$t(F/M)^{1/2}K^{-\varepsilon}\ll k$ can be assumed. In particular,
if $k\ll K^\varepsilon$, then $F\ll K^\varepsilon M/t^2 \ll K^\varepsilon$,
under Convention 1. We are led to a situation analogous to $(7.21)$, and invoking
the bound $(1.32)$ for $H_{j,k}\left({1\over2}+it\right)$, which is yet to be
proved, as well as $(5.33)$ or any uniform bound for $\tau_{j,k}(n)$, we can
settle this case. Note that Good's bound mentioned immediately after $(1.17)$ 
should not be utilised here, because it appears not to be uniform in the relevant
cusp forms, unlike Meurman's bound applied to $(7.21)$. At any event, we may
assume that
$k\gg K^\varepsilon$. Then, shifting the contour far to the left again, we are led
to the truncation
$$
\max\left\{K^\varepsilon,\,t(F/M)^{1/2}K^{-\varepsilon}\right\}\ll k\ll
t(F/M)^{1/2}K^\varepsilon.
\eqno(7.24)
$$
With this, we insert $(7.23)$ into $(7.22)$, and take the $f$-sum innermost. 
We see that $(7.22)$ is
$$
\eqalignno{
\ll GK^{1/2+\varepsilon}\sup_{U,F}&{1\over U}
\int_{G^2/K^{\varepsilon}}^\infty
\sum_{U\le k\le2U}\sum_{j=1}^{\vartheta(k)}\alpha_{j,k}
\left|H_{j,k}\left(\txt{1\over2}
+it_1\right)\right|^2\cr
&\times\left|\sum_{f=1}^\infty\phi_2(f)\phi(fx)
{d(f)\tau_{j,k}(f)\over f^{{1/2}+i(t_2-t_1)}}
\right|{dx\over x^{5/4}}, &(7.25)
}
$$
where the lower bound for $x$ comes from $(7.16)$.
To this we apply $(1.32)$--$(1.33)$ and Lemma 8. More precisely, 
to one factor of $\left|H_{j,k}\left({1\over2}+it_1\right)\right|^2$ we apply
$(1.32)$, which is possible, because
$t\gg U^{3/2}$ with $U$ in the range $(7.24)$; and
to the remaining part of $(7.25)$ we apply $(1.33)$ and Lemma 8. We find that
$(7.25)$ is
$$
\ll K^{1+\varepsilon}t^{{1/3}}\sup_{U,F}{1\over
U}\left\{\left(U^2+t^{2/3}\right)
\left(U^2+F\right)\right\}^{1/2}\ll
\left(Kt^{2/3}+t^{4/3}\right)^{1+\varepsilon}.
\eqno(7.26)
$$
This finishes the treatment of the contribution of $\e{H}_2$ to $(7.7)$, up to
the proof of $(1.32)$--$(1.33)$. 
\medskip
Now, we move to the contribution to $(7.7)$ of $\e{H}_4$ defined by $(2.13)$. 
This is equal to
$$
\sum_{n=1}^\infty {\phi(n)\sigma_{2it_1}(n)\over
n^{{1/2}+it_2}}\sum_{m=1}^{n-1} m^{-{1/2}}d(m)d(n-m)\Psi^-(m/n;h).
\eqno(7.27)
$$
We should remark first that the present choice $(7.2)$ of the function $h$
allows us to move the vertical line of $(2.8)$ to the left freely as far as
$a>-P_0-{1\over2}$ and $a\ne -{1\over2},\ldots, -P_0+{1\over2}$. This can be
seen by a simple extension of the argument on [23, p.\ 113 and p.\ 121].
Also we may replace $\Psi^-(x;h)$ by
$$
\eqalignno{
GK\int_0^\infty&\left\{\int_{(a)}x^s(y(y+1))^{s-1}
{\Gamma^2({1\over2}-s)\over \Gamma(1-2s)\cos(\pi s)}ds\right\}\cr
&\times\left({y\over y+1}\right)^{i\delta K}
\exp\left(-\left({G\over2}\log{y\over y+1}\right)^2\right)dy,&(7.28)
}
$$
with $\delta =\pm1$, and further by
$$
\eqalignno{
GK\int_{G/\log K}^\infty&\left\{\int_{a-i(\log K)^2}^{a+i(\log K)^2}
{(xy^2)^s\Gamma^2({1\over2}-s)\over \Gamma(1-2s)\cos(\pi s)}ds\right\}
\exp\left(-\delta iK/y-\txt{1\over 4}(G/y)^2\right){dy\over y^2}.\quad\qquad
&(7.29)
}
$$
In particular,
$$
\Psi^-(x;h)\ll K(xG^2)^{-P_0}.\eqno(7.30)
$$
Hence, taking $P_0$ sufficiently large, we may consider, instead of $(7.27)$, 
$$
\sum_{n=1}^\infty {\phi(n)\sigma_{2it_1}(n)\over
n^{{1/2}+it_2}}\sum_{m\ll K^\varepsilon M/G^2 }
m^{-{1/2}}d(m)d(n-m)\Psi^-(m/n;h).
\eqno(7.31)
$$
We may further replace this by
$$
\sum_{f=1}^\infty {\phi_2(f)d(f)\over
f^{1+it_2}}\sum_{n=1}^\infty d(n)\sigma_{2it_1}(n+f)V_f^-(n/f),
\eqno(7.32)
$$
where $\phi_2$ is as in $(7.14)$, with $(7.16)$, and
$$
V_f^{-}(x)=\phi(f(x+1))(x+1)^{-{1/2}-it_2}
\Psi^{-}\left((x+1)^{-1};h\right).\eqno(7.33)
$$
\medskip
Now, we shall proceed with
$$
a=0.\eqno(7.34)
$$
We apply $(2.17)$ to $(7.32)$, with $\alpha=0,\, \beta=2it_1$. We consider
first  the contribution of $D_r$. We need to take a limit  on the
right side $(2.22)$, but obviously this procedure can be ignored. Then, for
instance, the leading term of $(2.22)$ yields the expression
$$
\eqalignno{
&GK\int_{-i(\log K)^2}^{i(\log K)^2}
{\Gamma^2({1\over2}-s)\over\Gamma(1-2s)\cos(\pi s)}R_\delta(s;G,K)\cr
&\times \sum_{f=1}^\infty \phi_2(f)
{d(f)\sigma_{1+2it_1}(f)\over f^{1+it_2}}
\int_{G^2/K^{\varepsilon}}^\infty{ \phi(f(x+1))\over
(x+1)^{s+{1/2}-2it_1 +it_2}}dx ds, &(7.35)
}
$$
where
$$
R_\delta(s;G,K)=\int_{G/\log K}^\infty y^{2(s-1)}\exp\left(-\delta iK/y
-\txt{1\over4}(G/y)^2\right)dy.\eqno(7.36)
$$
Lemma 6 implies that the innermost integral of $(7.35)$ is negligibly small, and
the same assertion holds for the contribution of $D_r$ to $(7.32)$.
\medskip
Let us consider the contribution of $D_d$. 
We need to study the $\Psi_\pm$-functions defined by $(2.24)$--$(2.25)$,
with our current specifications. 
In view of $(7.29)$ with $a=0$ and $(7.33)$, we may deal instead with the
expression
$$
GK\int_{-i(\log K)^2}^{i(\log K)^2}
{\Gamma^2({1\over2}-s)\over \Gamma(1-2s)\cos(\pi s)}R_\delta(s;G,K)
\Delta_\pm(\kappa,f,s)ds,\eqno(7.37)
$$
respectively, where
$$
\eqalignno{
\Delta_+&(\kappa,f,s)={1\over4\pi i}
\int_{(\varepsilon)}
\cos(\pi s_1)\Gamma(s_1+i\kappa)\Gamma(s_1-i\kappa)
\Gamma^2\left(\txt{1\over2}-it_1-s_1\right)\cr
&\times 
\int_{G^2/K^{\varepsilon}}^\infty \phi(f(x+1))
x^{s_1+it_1-{1/2}}(x+1)^{-s-{1/2}-it_2}dxds_1,&(7.38)
 }
$$
$$
\eqalignno{
&\Delta_-(\kappa,f,s)={1\over4\pi i}\cosh(\pi\kappa)
\int_{(\varepsilon)}\sin\left(\pi(s_1+it_1)\right)
\Gamma(s_1+i\kappa)\Gamma(s_1-i\kappa)\cr
&\times \Gamma^2\left(\txt{1\over2}-it_1-s_1\right)
\int_{G^2/K^{\varepsilon}}^\infty \phi(f(x+1))
x^{s_1+it_1-{1/2}}(x+1)^{-s-{1/2}-it_2}dxds_1.&(7.39)
 }
$$
We observe that we may suppose, in both expressions, that
$|\kappa|\ll K^{1/\varepsilon}$. To see this, it suffices to adopt the
reasoning following $(7.19)$. Note that we have $s\ll \log^2K$.
\par 
Then, let us consider $\Delta_+$. The $\Gamma$-factor is
$$
\ll 
\exp\left(-\pi\max(|\lambda_1|,|\kappa|)-\pi|\lambda_1+t_1|\right),\quad \Im
(s_1)=\lambda_1.\eqno(7.40)
$$
Thus, the case $|\lambda_1+t_1|\ge{1\over2}t$ can readily be ignored. Otherwise,
the inner integral is obviously negligibly small by Lemma 6. 
Hence, $\Delta_+$ can be discarded.
\par
We turn to $\Delta_-$. We note that the factor
$\cosh(\pi\kappa)\sin\left(\pi(s_1+it_1)\right)$ is cancelled out by the
$\Gamma$-factor. Before shifting the $s_1$-contour, we shall show that we may
truncate it to $\lambda_1\ll K^\varepsilon$. In fact,  concerning the last
$x$-integral, we have
$$
{d\over dx}\left((\lambda_1+t_1)\log x-(\lambda+t_2)\log(x+1)\right)
={1\over x}\left(\lambda_1+{t\over
x+1}+O(K^\varepsilon)\right),\eqno(7.41)
$$
with $\Im(s)=\lambda$, since $\lambda\ll\log^2K$ and $|t_\nu-t|\ll
K^\varepsilon$ $(\nu=1,2)$. Provided $|\lambda_1|\gg K^\varepsilon$ and under
Convention 1, the absolute value of the right side of $(7.41)$ is $\gg|\lambda_1|/x$
because of the bound $t/(x+1)\ll K^\varepsilon$ which is implied by $(7.1)$
and $(7.3)$.  Thus Lemma 6 works, and the part of $\Delta_-$ with $|\lambda_1|\gg
K^\varepsilon$ can be discarded as claimed. With this, we shift the
$s_1$-contour far to the right, without encountering any pole, and obtain the
truncation $t(F/M)^{1\over2}K^{-\varepsilon}\ll |\kappa|$. Thus, if $\kappa\ll
K^\varepsilon$, then $F\ll K^\varepsilon$. That is, this case is settled by
Meurman's bound as before. On the other hand, if $|\kappa|\gg K^\varepsilon$, then
we may shift the $s_1$-contour far to the left under Convention 1, again without
encountering any pole, and come to the following analogue of
$(7.24)$:
$$
\max\left\{K^\varepsilon,\,t(F/M)^{1/2}K^{-\varepsilon}\right\}
\ll |\kappa|\ll t(F/M)^{1/2}K^\varepsilon.\eqno(7.42)
$$
Hence, the contribution of
$D_d$ to $(7.32)$ is
$$
\eqalignno{
\ll K^{1+\varepsilon}\sup_{U,F}{1\over U}&
\int_{G^2/K^{\varepsilon}}^\infty
\sum_{U\le \kappa_j\le2U}\alpha_j
\left|H_j\left(\txt{1\over2}
+it_1\right)\right|^2\cr
&\times\left|\sum_{f=1}^\infty\phi_2(f)\phi(f(x+1))
{d(f)\tau_j(f)\over f^{{1/2}+i(t_2-t_1)}}
\right|{dx\over x}, &(7.43)
}
$$
with $U$ in the range $(7.42)$.
The rest is similar to the discussion of $(7.25)$. This time we appeal instead 
to $(1.18)$--$(1.19)$ and Lemma 7. We find that $(7.43)$ is $\ll
\left(Kt^{2/3}+t^{4/3}\right)^{1+\varepsilon}$.
\medskip
As to the $D_c$ part, we follow the above reasoning, and get,
instead of $(7.43)$, the expression
$$
\eqalignno{
 K^{1+\varepsilon}\sup_{U,F}{1\over U}&
\int_{G^2/K^{\varepsilon}}^\infty
\int_{U}^{2U}{\left|\zeta\left(\txt{1\over2}+i(t_1+\kappa)\right)\right|^4
\over|\zeta(1+2i\kappa)|^2}\cr
&\times\left|\sum_{f=1}^\infty\phi_2(f)\phi(f(x+1))
{d(f)\sigma_{2i\kappa}(f)\over f^{{1/2}+i\kappa+i(t_2-t_1)}}
\right|d\kappa{dx\over x}. &(7.44)
}
$$
We have 
$\zeta\left({1\over2}+i(t_1+\kappa)\right)\ll t^{1/6}$ because of the upper
bound in $(7.42)$.  We apply the Mellin inversion of $\phi$ and $\phi_2$, 
with an appropriate truncation of the resulting new double integral. Then, the
sum over $f$ is expressed in terms of the zeta-function.
We may shift the two contours to the imaginary
axis without encountering any pole, because of
the lower bound in $(7.42)$. The inner integral of $(7.44)$  is
$\ll t^{2/3}UK^\varepsilon$. Hence $(7.44)$ itself is $\ll
K^{1+\varepsilon}t^{2/3}$.
\par
The treatment of the $D_h$ part is
analogous to that pertaining to $\Delta_+$, and the contribution is negligibly
small. 
\medskip
We conclude that the $\e{H}_2$ and $\e{H}_4$ parts
of $(7.7)$ are both $\ll\left(Kt^{2/3}+t^{4/3}\right)^{1+\varepsilon}$. 
Combined with the assertion adjacent to $(7.9)$, 
this ends the proof of Theorem 2
and  thus of Theorem 1, leaving one essential step yet to be confirmed. What 
remains is to prove $(1.33)$.
That is to be done in the next section.
\vskip 1cm
\centerline{\bf 8. Discussion}
\bigskip
\noindent
Here we shall first develop a brief proof of $(1.33)$, and that of Theorem 3. Then,
we shall observe a structure that makes such
extensions possible. With this,  the feasibility of further extensions will
be discussed.
\medskip
To prove $(1.33)$, we follow closely the argument of [13].  We first
look into the case $t\gg K^{3+\varepsilon}$ with a large $K$. The main
difference with the corresponding part of [13] is in that instead of a Voronoi
summation formula involving the coefficients $\tau_j(n)$ we work with its
counterpart for $\tau_{j,k}(n)$. This means replacing the Bessel
function $J_{2ir}(x)$ with a real $r\approx K$ by $J_{2k-1}(x)$ with an integer 
$k \approx K$. The argument in [13] relies on the fact that if 
$x \gg K^{2+\varepsilon}$, then
$$
{J_{2ir}(x)-J_{-2ir}(x)\over 2\sinh\pi r}\sim i\sqrt{{2\over\pi x}}
\sin\left(x-\txt{1\over4}\pi \right)\eqno(8.1)
$$
(see $(4.6)$). With the same assumption, we have
$$
J_{2k-1}(x) \sim (-1)^{k-1}\sqrt{{2\over\pi x}}
\sin\left (x-\txt{1\over4}\pi\right).
\eqno(8.2)
$$
if $k$ is a natural number (see $(8.5)$ below). The analogy is perfect as far
as the Bessel functions are concerned. Also,
to resulting sums involving the coefficients $\tau_{j,k}(n)$
we apply Lemma 8 in place of Lemma 7 that is used in [13]. Hence, the argument of
[13] can be repeated word by word if $t\gg K^{3+\varepsilon}$. 
\par
Thus we assume that $t\ll K^{3+\varepsilon}$ as well as
$K^\varepsilon\ll G\ll K^{1-\varepsilon}$. 
Then, we have, analogously to $(7.5)$, that
$$
H_{j,k}\left(\txt{1\over2}+it\right)
\ll (\log K)^2\sum_{M\ll K+t}
\int_{\gamma^{-1}-i\gamma^2}^{\gamma^{-1}+i\gamma^2}
\left|\sum_{m=1}^\infty
\phi(m;M)\tau_{j,k}(m)m^{-{1/2}-it-\xi}\right|
|d\xi|,\eqno(8.3)
$$
where $M$ runs over dyadic numbers.
With this, the case $t\ll K^{1+\varepsilon }$ is readily settled by Lemma 8. Hence,
it remains to consider the intermediate range
$K^{1+\varepsilon } \ll t \ll K^{3+\varepsilon }$. Here the proof of Lemma 8 is relevant. Thus, we are to deal with
$$
\sum_{\ell}{1\over\ell}\sum_{m,n}{\phi_0(m)
\overline{\phi_0(n)}\over\sqrt{mn}}\left(m/
n\right)^{it} S(m,n;\ell)(h_1)^\circ\left({4\pi}\sqrt{mn}/\ell\right),\eqno(8.4)
$$
where $\phi_0$ is as in $(7.6)$, $(h_1)^\circ$ the expression $(3.34)$, and the
truncation $(3.12)$ has already been applied, but with $N$ being replaced by $M\ll
t$. The Kloosterman sums are expanded according to their definition, and the
assertion $(3.36)$ is invoked. Then we end up with a double exponential sum over
$m$ and $n$,  essentially the same as the corresponding sum in [13]. This
ends the proof of
$(1.33)$. Consequentially, we have finished the proof of Theorem 1.
\medskip
As to the proof of Theorem 3, it depends solely on the observation that the
procedure developed in Section 4 is as a matter of fact a reduction
of the original problem to additive divisor sums. 
Applied to the left side of $(1.34)$, this argument leads us
to exactly the same additive divisor sums, albeit there exist differences coming
from the use of Lemma 2 in place of Lemma 1 and from that $T\approx (K+t)^2$.
There is virtually no difference in terms of asymptotics. This is
endorsed by the truncation $(3.12)$, which is applicable to the present situation 
as mentioned above, and by the formula
$$
J_{2k-1}(x)\sim (-1)^{k-1}\sqrt{2\over \pi
x}\sin\left(\omega(ik,x)-\txt{1\over4}\pi\right).\eqno(8.5)
$$
The former fact corresponds to $(4.3)$, and the latter to $(4.6)$, respectively.
The rest of the proof is the same as that of Theorem 1. In fact, it is slightly
simpler, because the second condition in $(6.1)$ is unnecessary, due to the fact
that we have $ik$ in $(8.5)$ in place of $r$ in $(4.6)$.
\bigskip
We shall expand our observation about the r\^ole of additive
divisor sums. To this end, we return to $(3.40)$. The $L$-series that yields the
Dirichlet series on the right is associated with 
the Rankin--Selberg convolution of the Eisenstein series and the relevant cusp
form. The divisor function there is  a Fourier coefficient of an
automorphic function. The structure of our subsequent reasoning, which is
admittedly involved, could be summarised as follows:
\medskip
\item{(1)} Appearance at $(4.2)$ of Kloosterman sums via Lemma 1 
\item{(2)} Basic truncation $(4.3)$ of moduli of Kloosterman sums 
\item{(3)} Application of the Vorono{\"\i} sum formula at $(4.11)$
\item{(4)} Another basic truncation at $(4.21)$
\item{(5)} Appearance  at $(4.27)$ of additive divisor sums
\item{(6)} Application at $(5.1)$ of the spectral decomposition $(2.27)$
\item{(7)} Truncation of the spectral range at $(5.20)$/$(6.24)$
\item{(8)} Appearance at $(6.33)$ of a simpler spectral sum
\item{(9)} Reduction at $(6.37)$ to a hybrid moment 
\medskip
\noindent
Note that Step $(3)$ is performed upon the divisor function that is
never of our main concern at $(3.40)$. The subsequent analysis is, however,
wholly relevant to these Fourier coefficients of the Eisenstein series. It is
true that Kloosterman sums replace Fourier coefficients of original cusp forms and
thus the latter objects are actually playing a r\^ole in the background.
Nevertheless, those operations following $(3)$ are made possible because of
the presence of the divisor function. Moreover, the decisive step $(8)$ is due
solely to $(5)$. In other words, the divisor function is indeed the protagonist of
our scenario, despite its  obscure entrance at $(3.40)$. Or perhaps more
correctly, an orchestration of automorphic waves conducted by the sum formulas due
to Bruggeman, Kuznetsov, and Petersson  makes it possible for the divisor function
to conjure the uniform subconvexity bounds $(1.1)$ and $(1.2)$.
\medskip
Now, if $(3.40)$ can be regarded as a statement concerning a Rankin--Selberg
convolution, then what has been developed above could be a typical instance of a
general mechanism arising from automorphy; by no means a
serendipity. We shall
indicate, with a plausible inference, that this should be the case.
\medskip
Thus, let $\psi$ be a Hecke invariant cusp form, either holomorphic or real
analytic. Let $\tau_\psi(n)$ be its Hecke
eigenvalue. We are interested in bounding the Rankin--Selberg $L$-function
$$
L(s,\psi\otimes\psi_j)=\zeta(2s)\sum_{n=1}^\infty
\tau_\psi(n)\tau_j(n)n^{-s}\eqno(8.6)
$$
on the critical line. Note that the function $(s-1)L(s,\psi\otimes\psi_j)$ is
entire, and also that one may naturally replace $\psi_j$ by $\psi_{j,k}$, and
proceed analogously. 
\par
We need to treat the expression
$$
\sum_{K\le\kappa_j\le K+G}\alpha_j
\left|\sum_{n=1}^\infty \phi(n)\tau_\psi(n)\tau_j(n) n^{-{1/2}-it}\right|^2,
\eqno(8.7)
$$
where $(1.26)$ is effective, and $\phi$ as in $(4.2)$ with $M\ll T_\psi$, where
 $T_\psi$ is defined analogously to $(3.38)$.
We may apply Steps $(1)$ and $(2)$ without any change. The third step is
equivalent to an appeal to the functional equation for the Hecke--Estermann
zeta-unction
$$
\sum_{n=1}^\infty\tau_\psi(n)\exp\left(2\pi iqn/\ell\right)n^{-s},
\quad (q,\ell)=1,\eqno (8.8)
$$
which is an extension of $(4.14)$, and a consequence of the automorphy of $\psi$.
Essentially the same as $(4.17)$ comes out, with $d$ being replaced by $\tau_\psi$.
Here might, however, arise a problem relevant to the change in the function $I$,
which should be taken into account if the uniformity in $\psi$ is to be
maintained. The same can be said about the extension of Step $(4)$. Step $(5)$ is
now with the sum
$$
\sum_{n=1}^\infty \tau_\psi(n)\tau_\psi(n+f)W(n/f).\eqno(8.9)
$$
When $\psi$ is holomorphic, there exists a complete analogue of 
Lemma 5 that is due to the second named author (implicit in [22]). Hence this case
should not cause any extra difficulty as far as Step $(8)$. With a real
analytic $\psi$, there might arise a new issue, because we lack any complete
extension of Lemma 5 to this case. There exists, however, a relevant result, an
asymptotic extension due to the first named author [9]. That might serve well for
our purpose. Despite this, we should better try to achieve a complete extension of
Lemma 5 to the real analytic case, mainly for the sake of a fuller understanding
of this fascinating mechanism.  In fact, such a programme is being undertaken by
the second named author (see [24]); the key seems to be the harmonic
analysis on the Lie group $\r{PSL}_2(\B{R})$. Thus, we may envisage with a good
reason that we could go through Steps $(5)$--$(7)$ in the new context as well; 
that is, an analogue of $(6.33)$ should hold with $\psi$ in general. 
There the factor $H_j^2\!\left({1\over2}\right)$ is 
to be replaced by the inner product
$\langle |\psi|^2,\psi_j\rangle$ or a quantity closely related to, with an
appropriate normalisation of the metric. We need an analogue for $\langle
|\psi|^2,\psi_j\rangle$  of the spectral fourth moment of
$H_j\!\left({1\over2}\right)$. Such a result,
in fact the spectral mean square of the inner product, is proved by Good [4].
Also, its extension to the real analytic case is obtained by the first named
author [9]. Other parts of $(6.33)$ do not need to be changed substantially.
Therefore, it is highly probable that a counterpart of $(1.27)$, 
and consequentially a subconvexity bound 
$$
L\left(\txt{1\over2}+it,\psi\otimes\psi_j\right)\ll \kappa_j^{{2/3}+\varepsilon}
\eqno(8.10)
$$ 
be within our reach, at least when $t$ is relatively small
compared with $\kappa_j$. Indeed, we have proved already that this is the case when
$\psi$ is holomorphic, with a meaningful uniformity in $\psi$ and $t$. 
The situation with a real analytic $\psi$ should be analogous, 
though we have not worked out the details as yet. 
\par
It remains to ponder about a fuller analogue of Theorem 1. Here we are, however,
to realise that we were in a fortuitous situation with ${\cal S}(G,K)$.  A reason
why the hybrid mean value worked fine with $(6.33)$ is in that
the latter has the factor $H_j^2\!\left({1\over2}\right)$, as this fact was
exploited  to reach $(6.37)$. Such a splitting of the corresponding factor 
$\langle |\psi|^2,\psi_j\rangle$ does not appear to be
possible in general. Thus, we should better stop our plausible inference here.
It should, however, be added that there are other directions of the extension. For
instance, we may replace the group $\varGamma$ by $\varGamma_0(q)$; and the
twist of $H_j(s)$ with a Dirichlet character can be treated with the same strategy
as above, taking into account the uniformity in the modulus of the character.
Another possibility is to include the Bianchi groups. All basic machineries
needed for this purpose are laid out in [2]. 
\medskip
Finally, we stress that there exists a possibility that one might
come to Step $(8)$ directly from the original spectral sum. That is, the use of
Kloosterman sums and the Voronoi sum formula could be avoided altogether. This is
suggested by the recent work [3], where the spectral decomposition of the fourth
power moment of the Riemann zeta-function is grasped as a special instance of the
same of a Poincar\'e series on the group $\r{PSL}_2(\B{R})$, yielding a new
approach to the subject closely related to ours. In the perspective thus opened,
the functional equations and the Bessel transforms which are in the core of our
analysis developed above are understood to be realisations of the action of the
Weyl element of the group under various circumstances.
To this and the above observation on Rankin--Selberg $L$-functions
we shall return elsewhere.
\medskip
\noindent
{\bf Concluding remark.} After finishing the 
present work, we found that Sarnak had developed in [25] an approach 
to the subconvexity bound of Rankin--Selberg $L$-functions. He worked
mainly with holomorphic cusp forms; nevertheless, the initial stage of his approach
is analogous to ours in the sense that the corresponding steps up to $(5)$
are observable there, though with a different configuration. 
It is indeed hard to conceive any other way to take. However, from the stage
corresponding to $(6)$ on, Sarnak's strategy differs considerably from ours, and 
the bound that his method gives rise to is tangibly weaker than our assertion
pertaining to $(8.10)$, as far as the full modular group is concerned. A talk on
this subject and in fact a summary of the present article were delivered by the
authors at the Tagung `Theory of the Riemann Zeta and Allied Functions'
(Mathematisches  Forschungsinstitut Oberwolfach, September 20, 2004).
\vskip 1cm
\centerline{\bf References}
\bigskip
\item{[1]} R.W.\ Bruggeman. Fourier coefficients of
cusp forms. Invent.\ math., {\bf 45} (1978), 1--18.
\item{[2]} R.W. Bruggeman and Y. Motohashi. Sum formula  for Kloosterman 
sums and the  fourth moment of the Dedekind zeta-function   over  the Gaussian
number field. Functiones et Approximatio, {\bf 31} (2003),  7--76. 
\item{[3]} ---. A new approach to the
spectral theory of the fourth moment of the Riemann zeta-function. 
To appear in Crelle's J.
\item {[4]} A.\ Good. The square mean of Dirichlet series associated with
cusp forms.  Mathematika, {\bf 29} (1982), 278-295.
\item{[5]} A. Ivi\'c. On sums of Hecke series in short intervals.
J. Th\'eorie des Nombres de Bordeaux, {\bf 13} (2001), 554--568.
\item{[6]} H. Iwaniec. Fourier coefficients of cusp forms and the Riemann
zeta-function. Exp.\ No.\thinspace18, S\'em.\ Th.\ Nombres, Univ.\ Bordeaux 
1979/80.
\item{[7]} ---  Small eigenvalues of Laplacian for
$\varGamma _0(N)$.  Acta Arith., {\bf 56} (1990), 65--82.
\item{[8]} M. Jutila. A Method in the Theory of Exponential Sums. 
Tata IFR, Lect.\ in Math.\ Phy., {\bf 80}, Springer Verlag 1987.
\item{[9]} ---  The additive divisor problem and its analogs
for Fourier coefficients of cusp forms. I. Math. Z., {\bf 223} (1996), 435--461;
II. ibid {\bf 225} (1997), 625-637.
\item{[10]} --- Mean values of Dirichlet series via Laplace
transforms.  In: Analytic Number Theory, Proc.\ 39th Taniguchi Intern.\ Symp.\
Math., Y. Motohashi (Ed.), Cambridge Univ.\ Press 1997, pp.\ 169--207.
\item{[11]} --- On spectral large sieve inequalities. Functiones et Approximatio,
{\bf 28}  (2000), 7--18.
\item{[12]} --- The fourth moment of central values of Hecke
series. In: Number Theory, Proc.\ Turku Symposium on Number
Theory in Memory of Kustaa Inkeri, M. Jutila and T.
Mets\" ankyl\" a (Eds.), Walter de Gruyter 2001,
pp.\ 167--177.
\item{[13]} --- The spectral mean square of Hecke $L$-functions on the 
critical line. To appear in Publ.\ de l'Institut Math.\ Belgrade.
\item{[14]}  M. Jutila and Y. Motohashi. A note on the mean value
of the zeta and $L$-functions. XI. Proc. Japan Acad., {\bf 78A} (2002), 1-6.
\item{[15]} S. Katok and P. Sarnak. Heegner points, cycles and
Maass forms. Israel J. Math., {\bf 84} (1993), 193--227.
\item{[16]} N.V.\ Kuznetsov.\enskip Petersson hypothesis for forms of weight zero
and Linnik  hypothesis. Khabarovsk Complex Res.\ Inst.\
Acad.\ Sci.\ USSR. Preprint 1977. (Russian)
\item{[17]} --- Convolution of the Fourier coefficients of 
the Eisenstein--Maass series. Zap.\ Nauchn.\ Sem.\ LOMI, {\bf
129} (1983),  43--84. (Russian)
\item{[18]}  N.N. Lebedev. Special Functions and their Applications.
Dover 1972.
\item{[19]} T. Meurman. On the order of the Maass $L$-function on the
critical line. Coll.\ Math.\ Soc.\ J\'anos Bolyai, {\bf 51} (1989),
325--354.
\item{[20]} Y. Motohashi. An explicit formula for the
fourth   power mean of the Riemann zeta-function. Acta Math., {\bf
170} (1993),   181--220.
\item{[21]} --- The binary additive divisor problem. 
Ann.\ Sci.\ \'Ecole Norm.\ Sup., (4) {\bf27} (1994), 529--572.
\item{[22]} --- The mean square of Hecke
$L$-series attached to holomorphic cusp forms. RIMS Kyoto Univ.
Kokyuroku, {\bf 886} (1994), 214--227.
\item{[23]} --- Spectral theory of the Riemann
zeta-function. Cambridge Tracts in Math.,\ {\bf 127}, Cambridge
Univ.\ Press 1997.
\item{[24]} --- A note on the mean value of the zeta and
$L$-functions. XIV. Proc.\ Japan Acad., {\bf 80A} (2004), 28--33. 
\item{[25]} P. Sarnak. Estimation of Rankin--Selberg $L$-functions
and quantum unique ergodicity. J. Funct.\ Analy., {\bf 184} (2001),
419--453.
\item{[26]}  G.N. Watson.  A treatise on the theory of Bessel
functions. Cambridge Univ.\ Press 1996.

\bigskip
\noindent
\font\small=cmr8 {\small Matti Jutila
\par\noindent Department of Mathematics, University of Turku,
\par\noindent FIN-20014 Turku, Finland
\par\noindent Email: jutila@utu.fi}
\medskip
\noindent {\small Yoichi Motohashi
\par\noindent Department of Mathematics, College of Science and
Technology, Nihon University,
\par\noindent Surugadai, Tokyo 101-8308, Japan
\par\noindent Email: ymoto@math.cst.nihon-u.ac.jp,
am8y-mths@asahi-net.or.jp }

\bye